\documentclass[a4paper,12pt,nosumlimits]{amsart}
\usepackage{latexsym,amsmath,amssymb,color,enumerate}
\usepackage[english]{babel}
\numberwithin{equation}{section}
\allowdisplaybreaks[4]

\usepackage{color}

\title[Calculi, Hodges and Laplacians]
{Calculi, Hodge operators  and Laplacians \\ on a quantum Hopf fibration}
\date{v3:  3 Oct. 2012}
\author{Giovanni Landi}
\address{Dipartimento di Matematica e Informatica, Universit{\`a} di Trieste, 
Via  A. Valerio 12/1, I-34127, Trieste, Italy;  
and INFN, Sezione di Trieste, Trieste, Italy.} 
\author{Alessandro Zampini}
\address{Mathematisches Institut der Ludwig - Maximilians - Universit\"at, Theresienstra\ss e 39, 80333 M\"unchen, Germany } 
\email{landi@univ.trieste.it, zampini@math.lmu.de} 

\newtheorem{theo}{Theorem}[section]
\newtheorem{lemm}[theo]{Lemma}
\newtheorem{prop}[theo]{Proposition}

\newtheorem{rema}[theo]{Remark}

\newcommand{\ii}{\mathrm{i}}
\newcommand{\nn}{\nonumber}

\newcommand{\dd}{{\rm d}}
\newcommand{\ca}{\mathcal{A}}

\newcommand{\ch}{\mathcal{H}}
\newcommand{\cl}{\mathcal{L}}
\newcommand{\cn}{\mathcal{N}}
\newcommand{\cp}{\mathcal{P}}
\newcommand{\cb}{\mathcal{B}}
\newcommand{\cq}{\mathcal{Q}}
\newcommand{\oca}[1]{\Omega^{#1}(\ca)}
\newcommand{\och}[1]{\Omega^{#1}(\ch)}
\newcommand{\cf}{\mathcal{F}}
\newcommand{\cR}{\mathcal{R}}

\newcommand{\cu}{\mathcal{U}}        
\newcommand{\SU}{\mathrm{SU}_q(2)}  
\newcommand{\ASU}{\ca(\mathrm{SU}_q(2))}  
\newcommand{\sq}{\mathrm{S}^2_{q}}  
\newcommand{\Asq}{\ca(\mathrm{S}^2_{q})}  
\newcommand{\su}{\cu_q(\mathfrak{su}(2))}  

\newcommand{\eps}{\varepsilon}      
\newcommand{\co}[2]{#1_{(#2)}}      
\newcommand{\hs}[2]{\left\langle #1,#2\right\rangle}  
\newcommand{\hp}[2]{\left\{ #1,#2\right\}} 

\newcommand{\ket}[1]{\left | #1 \right\rangle }
\newcommand{\bra}[1]{\left\langle #1 \right |}
\newcommand{\half}{{\tfrac{1}{2}}}
\newcommand{\lt}{{\triangleright}}    
\newcommand{\rt}{{\triangleleft}}
\newcommand{\IC}{{\mathbb C}} 
\newcommand{\IR}{{\mathbb R}} 
\newcommand{\IN}{{\mathbb N}} 
\newcommand{\IZ}{{\mathbb Z}} 
\DeclareMathOperator{\Ad}{Ad}       
\DeclareMathOperator{\id}{id}       
\DeclareMathOperator{\U}{U}       

\newcommand{\abs}[1]{\left|#1\right|}

\newcommand{\figureheight}{8cm}
\newcommand{\putfig}[2]{\begin{figure}[htp]
        \special{isoscale c:/itex/texfig/#1.wmf, \the\hsize \figureheight}
        \vspace{\figureheight}
        \caption{#2}\label{fig:#1}
        \end{figure}}
\newcommand{\pictureheight}{4cm}
\newcommand{\putpicture}[2]{\begin{figure}[htp]
        \special{isoscale c:/itex/texfig/#1.wmf, \the\hsize \pictureheight}
        \vspace{\pictureheight}
        \caption{#2}\label{fig:#1}
        \end{figure}}

\newcommand{\beqa}{\begin{eqnarray}}
\newcommand{\eeqa}{\end{eqnarray}}
\newcommand{\beq}{\begin{equation}}
\newcommand{\eeq}{\end{equation}}

\newcommand{\Psin}{\Psi^{\left(n\right)}}

\newcommand{\mn}{\abs{n}}
\newcommand{\mj}{\abs{j}}

\newcommand{\qpp}{\mathfrak{p}^{\left(n\right)}}

\newcommand{\an}{\mathrm{a}^{(n)}}
\newcommand{\An}{\mathrm{A^{(n)}}}
\newcommand{\Ans}{\mathrm{A}_{s}^{(n)}}

\newcommand{\qum}{q^{\frac{1}{2}}}
\newcommand{\qun}{q^{-\frac{1}{2}}}

\newcommand{\sigmat}{\tilde{\sigma}}

\begin{document}

\thispagestyle{empty}
\begin{abstract}
\noindent
We describe Laplacian operators on the quantum group $\SU$ equipped with  the four dimensional bicovariant differential calculus of Woronowicz as well as on the quantum homogeneous space $\sq$ with the restricted left covariant three dimensional differential calculus. This is done by giving a family of Hodge dualities on both the exterior algebras of $\SU$ and $\sq$. We also study gauged Laplacian operators acting on sections of line bundles over the quantum sphere.  
\end{abstract}

\maketitle

\tableofcontents

\linespread{1.2}

\section{Introduction}

We continue our program devoted to Laplacian operators on quantum spaces with the study of such operators on the quantum (standard) Podle\'s sphere $\sq$ and their coupling with gauge connections on the quantum principal $\U(1)$-fibration $\Asq\hookrightarrow\ASU$. While in \cite{lareza} one worked with a left $3D$ covariant differential calculus on $\SU$ and its restriction to the (unique) $2D$ left covariant differential calculus on the sphere $\sq$, in the present paper we use the somewhat more complicate $4D_{+}$ bicovariant calculus on $\SU$ introduced in \cite{wor89} and its restriction to a $3D$ left covariant calculus on the sphere $\sq$. 

Laplacian operators on all Podle\'s spheres, related to the $4D_{+}$ bicovariant calculus on $\SU$ were already studied in \cite{poddc}. Our contribution to Laplacian operators comes from the use of Hodge $\star$-operators on both the manifold of $\SU$ and $\sq$ that we introduce by improving and diversifying on existing  definitions. 

We then move on to line bundles on the standard sphere $\sq$ and to a class of operators on such bundles that are `gauged' with the use of a suitable class of connections on the principal bundle $\Asq\hookrightarrow\ASU$ and of the corresponding covariant derivatives on (module of sections of) the line bundles. These gauged Laplacians are completely diagonalized and are split in terms of a Laplacian operator on the total space $\SU$ of the bundle minus vertical operators, paralleling what happens on a classical principal bundle  (see e.g. \cite[Prop.~5.6]{bgv}) and on the Hopf fibration of the sphere $\sq$ with calculi coming from the left covariant one on $\SU$ as shown in \cite{lareza,ale09}

In \S\ref{s:qsb} we describe all we need of the principal fibration $\Asq\hookrightarrow\ASU$ and associated line bundles over $\sq$. We also give a systematic description of the differential calculi we are interested in, the
$4D_{+}$ bicovariant calculus on $\SU$ and its restriction to a $3D$ left covariant calculus on the sphere $\sq$.
A thoughtful construction of Hodge $\star$-dualities on $\SU$ is in \S\ref{s:Hop}, while that on $\sq$ is in 
\S\ref{s:hl3}. These are used in \S\ref{se:L} for the definition of Laplacian operators. A digression on connections on the principal bundle and of covariant derivatives on the line bundles is in \S\ref{s:ccd} and the following  
\S\ref{GLoLB} is devoted to the corresponding gauged Laplacian operators on modules of sections of the line bundles. To make the paper relatively self-contained it concludes with two appendices, 
App.\ref{ass:a1} giving general facts on differential calculi on Hopf algebras and App.\ref{ap:qpb} concerning with general facts on quantum principal bundles endowed with connections.

We like to mention that besides the constructions in \cite{hec99} and \cite{kmt}, examples of Hodge operators on the exterior algebras of the quantum homogeneous $q$-Minkowski and $q$-Euclidean spaces -- satisfying a covariance requirement  with respect to the action of the quantum groups $\mathrm{SO}_{q}(3,1)$ and $\mathrm{SO}_{q}(4)$ -- have been given in \cite{um94, majqe} using the formalism of braided geometry and with a construction of a $q$-epsilon tensor. On  the exterior algebra over the quantum planes $\IR_{q}^{N}$ a Hodge operator has been studied in \cite{gf}. The one on $\IR_{q}^{4}$ reduces to a Hodge operator on $\SU$ for both its two bicovariant calculi \cite{gfbis}.

\subsection*{Conventions and notations} When writing about connections and covariant derivatives we shall pay attention in keeping the two notions distinct: a connection will be a projection on a principal bundle while a covariant derivative will be an operator on section, both objects fulfilling suitable properties. For $q \neq 1$ the `$q$-number' is defined as
\begin{equation}
[x] = [x]_q := \frac{q^x - q^{-x}}{q - q^{-1}} ,
\label{eq:q-integer}
\end{equation} 
for any $x \in \IR$. For a coproduct
$\Delta$ we use the Sweedler notation
$\Delta(x)=x_{(1)}\otimes x_{(2)}$, with implicit summation. This is iterated to $(\id\otimes\Delta)\circ\Delta(x) = (\Delta\otimes\id)\circ\Delta(x) = x_{(1)}\otimes x_{(2)}\otimes x_{(3)}$, and so on. 

\subsection*{Acknowledgments}
We are grateful to S. Albeverio, L.S. Cirio and I. Heckenberger for comments and suggestions.
AZ thanks P. Lucignano for his help with Maple. 
GL was partially supported by the Italian Project `Cofin08--Noncommutative Geometry, Quantum Groups and Applications'. AZ gratefully acknowledges the support of the Max-Planck-Institut f\"ur Mathematik in Bonn, the Hausdorff Zentrum f\"ur Mathematik der Universit\"at Bonn, the Stiftelsen Blanceflor Boncompagni-Ludovisi f\"odd Bildt (Stockholm), the I.H.E.S. (Bures sur Yvette, Paris).  

\section{Prelude: calculi and line bundles on quantum spheres}\label{s:qsb}

We introduce the manifolds of the quantum group $\SU$ and its quantum homogeneous space $\sq$ -- the standard Podle\'s sphere. The corresponding inclusion $\Asq\hookrightarrow\ASU$ of the corresponding coordinate algebras is a (topological) quantum principal bundle. Following App.~\S\ref{ass:a1} we then equip $\ASU$ with a 4-dimensional bicovariant calculus, whose restriction gives a 3-dimensional left covariant calculus on $\Asq$. 


\subsection{Spheres and bundles}
For the  quantum group
$\SU$ its polynomial algebra  $\ASU$ is the unital $*$-algebra generated by elements $a$ and $c$, with relations
\begin{align}
\label{derel}
& ac=qca\quad ac^*=qc^*a\quad cc^*=c^*c , \nn \\
& a^*a+c^*c=aa^*+q^{2}cc^*=1 .
\end{align}
In the limit $q \to1$ one recovers the commutative coordinate algebra on the group manifold 
$\mathrm{SU(2)}$. The algebra $\ASU$ can be completed to a $C^*$-algebra in a usual way by considering all its admissible representations and the supremum (universal) norm on them \cite{wor87}.
For the sake of the present paper this is not necessary since we are interested in   Laplacian operators on $\SU$ (and on its homogeneous space, the quantum sphere) and their spectra. Thus we only exhibit a vector space basis for $\ASU$ in \eqref{bsphi} below, giving an analogue of the classical Wigner $D$-functions for the  $\mathrm{SU}(2)$ group, i.e.  matrix elements of its unitary irreducible (co)-representations.  Also, without loss of generality, the deformation parameter
$q\in\IR$ will be restricted to the interval $0<q<1$, the map $q \to q^{-1}$ giving isomorphic algebras. 

If we use the matrix 
$$ U = 
\left(
\begin{array}{cc} a & -qc^* \\ c & a^*
\end{array}\right) , 
$$
whose being unitary is equivalent to relations \eqref{derel}, the Hopf algebra structure for $\ASU$ is given by coproduct, 
antipode and counit:
$$
\Delta\, U = U \otimes U ,  \qquad S(U) = U^* , \qquad \eps(U) = 1 , 
$$
that is $\Delta(a)= a \otimes a - q c^* \otimes c$, and $\Delta(c)= c \otimes a + a^* \otimes c$; 
$S(a)=a^*$ and $S(c)=-qc$; $\eps(a)=1$ and $\eps(c)=0$ and their $*$-conjugated relations.

\medskip
The quantum universal envelopping algebra $\su$ is the unital Hopf $*$-algebra
generated as an algebra by four elements $K^{\pm 1},E,F$ with $K K^{-1}=1=K^{-1}K$ and relations: 
\beq 
K^{\pm1}E=q^{\pm1}EK^{\pm1}, \qquad 
K^{\pm1}F=q^{\mp1}FK^{\pm1}, \qquad  
[E,F] =\frac{K^{2}-K^{-2}}{q-q^{-1}} . 
\label{relsu}
\eeq 
The $*$-structure is
$K^*=K, \,  E^*=F $,
and the Hopf algebra structure is provided  by coproduct
$$\Delta(K^{\pm1}) =K^{\pm1}\otimes K^{\pm1}, \quad
\Delta(E) =E\otimes K+K^{-1}\otimes E,  \quad 
\Delta(F)
=F\otimes K+K^{-1}\otimes F,$$
while the antipode is 
$S(K) =K^{-1}, \,
S(E) =-qE, \, 
S(F) =-q^{-1}F$ 
and the counit reads
$\varepsilon(K)=1, \,\varepsilon(E)=\varepsilon(F)=0$. The quadratic element 
\beq
C_{q}=\frac{qK^2-2+q^{-1}K^{-2}}{(q-q^{-1})^2}+FE-\tfrac{1}{4}
\label{cas}
\eeq
is a quantum Casimir operator that generates the centre of $\su$.

\medskip
The Hopf $*$-algebras  $\su$ and $\ASU$ are dually paired. The $*$-compatible bilinear mapping $\hs{~}{~}:\su\times\ASU\to\IC$ is on the generators given by
\begin{align}
&\langle K^{\pm1},a\rangle=q^{\mp 1/2}, \qquad  \langle K^{\pm1},a^*\rangle=q^{\mp 1/2},   \nn\\
&\langle E,c\rangle=1, \qquad \langle F,c^*\rangle=-q^{-1}, \label{ndp}
\end{align}
with all other couples of generators pairing to zero. This pairing is proved \cite{KS97} to be non-degenerate. The algebra $\su$ is recovered as a $*$-Hopf subalgebra in the dual algebra $\ASU^o$, the largest Hopf $*$-subalgebra contained in the vector space dual $\ASU^{\prime}$.

There are \cite{wor87} $*$-compatible canonical commuting actions of $\su$ on $\ASU$: 
$$
h \lt x := \co{x}{1} \,\hs{h}{\co{x}{2}}, \qquad
x \rt  h := \hs{h}{\co{x}{1}}\, \co{x}{2}. 
$$
On powers of generators one computes, for $s\in\,\IN$, that 
\beq
\label{lact}
\begin{array}{lll} 
K^{\pm1}\triangleright a^{s} =q^{\mp\frac{s}{2}}a^{s} & 
F\triangleright a^{s} =0 &
E\triangleright a^{s} =-q^{(3-s)/2} [s] a^{s-1} c^{*} \\
K^{\pm1}\triangleright a^{* s} =q^{\pm\frac{s}{2}}a^{* s} & 
F\triangleright
a^{*s} =q^{(1-s)/2} [s] c a^{* s-1} &
E\triangleright a^{* s} =0 \\
K^{\pm1}\triangleright c^{s} =q^{\mp\frac{s}{2}}c^{s} &
F\triangleright c^{s} =0 &
E\triangleright c^{s} =q^{(1-s)/2} [s]  c^{s-1} a^* \\
K^{\pm1}\triangleright c^{* s} =q^{\pm\frac{s}{2}}c^{* s} &
F\triangleright c^{*s} =-q^{-(1+s)/2} [s] a c^{*s-1}  &
E\triangleright c^{* s} =0;
\end{array}
\eeq
and:
\beq
\label{ract}
\begin{array}{lll}
 a^{s}\triangleleft K^{\pm1} =q^{\mp\frac{s}{2}}a^{s} & 
a^{s}\triangleleft F =q^{(s-1)/2} [s] c a^{s-1} &
a^{s}\triangleleft E =0 \\
 a^{* s}\triangleleft K^{\pm1} =q^{\pm\frac{s}{2}}a^{* s} & 
a^{* s}\triangleleft F =0 &
a^{*s}\triangleleft E =-q^{(3-s)/2} [s] c^{*}a^{*s-1} \\
c^{s}\triangleleft K^{\pm1} =q^{\pm\frac{s}{2}}c^{s}  & 
c^{s}\triangleleft F =0  &
c^{s}\triangleleft E =q^{(s-1)/2} [s] c^{s-1} a \\
c^{* s}\triangleleft K^{\pm1} =q^{\mp\frac{s}{2}}c^{* s} &
c^{* s}\triangleleft F =-q^{(s-3)/2} [s] a^{*}c^{*s-1} &
c^{* s}\triangleleft E =0. 
\end{array}
\eeq

\medskip
Consider the algebra $\ca(\U(1)):=\IC[z,z^*] \big/ \!\!<zz^* -1>$. The map  
\beq  \label{qprp}
\pi: \ASU \, \to\,\ca(\U(1)) , \qquad  
\pi\,\left(
\begin{array}{cc} a & -qc^* \\ c & a^*
\end{array}\right):=
\left(
\begin{array}{cc} z & 0 \\ 0 & z^*
\end{array}\right)
\eeq 
is a surjective Hopf $*$-algebra homomorphism. As a consequence, $\U(1)$
is a quantum subgroup of $\SU$ with right coaction:
\beq 
\delta_{R}:= (\id\otimes\pi) \circ \Delta \, : \, \ASU \,\to\,\ASU \otimes
\ca(\U(1)) . \label{cancoa} 
\eeq 
The coinvariant elements for this
coaction, 
elements $b\in\ASU$ for which $\delta_{R}(b)=b\otimes 1$,  
form the algebra of the standard Podle\'s sphere $\Asq\hookrightarrow\ASU$.  This inclusion gives a topological quantum principal bundle, following the formulation reviewed in appendix~\S\ref{ap:qpb}. 

\medskip
The above right $\U(1)$ coaction on $\SU$ is dual to the left action of the element $K$, and allows one    \cite{maetal} to give a decomposition 
\beq\label{dcmp}
\ASU=\oplus_{n\in\IZ} \cl_{n}
\eeq
in terms of  $\Asq$-bimodules defined by 
\beq
\label{libu} 
\cl_{n} := \{x \in \ASU ~:~ K \lt x = q^{n/2} x\quad\Leftrightarrow\quad\delta_{R}(x)=x\otimes z^{-n}\},
\eeq 
with $\Asq = \cl_{0}$.  It is easy to see (cf. \cite[Prop. 3.1]{mgw})  that  $\cl_{n}^{*} =  \cl_{-n}$ and $\cl_{n}\cl_{m} = \cl_{n+m}$. Also
\beq
E \lt \cl_{n} \subset \cl_{n+2}, \qquad 
F \lt \cl_{n} \subset \cl_{n-2}, \qquad
\cl_{n} \rt  u \subset \cl_{n},
\label{rellb}
\eeq
for any $u\in \su$. The bimodules $\cl_{n}$ will be described at length later on when we endow them with connections. 
Here we only mention that the bimodules $\cl_{n}$ have a vector space decomposition (cf. e.g. \cite{maj95}):
\beq
\label{decoln}
\cl_{n}:=\bigoplus_{J=\tfrac{|n|}{2}, \tfrac{|n|}{2} +1,
\tfrac{|n|}{2} +2, \cdots}V_{J}^{\left(n\right)}, 
\eeq 
where $V_{J}^{\left(n\right)}$ is  the spin
$J$ (with $J\in\half\IN$) irreducible $*$-re\-pre\-sen\-ta\-tion spaces for the right action of  $\su$, and basis  elements  
\beq
\phi_{n,J,l}= (c^{J-n/2} a^{*J+n/2})\rt E^l
\label{bsphi}
\eeq
with $n\in\IZ,\,J=\tfrac{|n|}{2}+\IN, \,l=0,\ldots,2J$.

\subsection{The 4D exterior algebra over the quantum group $\SU$}\label{se:4dc} 

Following the formulation reviewed in App.\ref{ass:a1}, we present here the exterior algebra over the  so called $4D_{+}$ bicovariant calculus on $\SU$, which was introduced  as a first order differential calculus  in \cite{wor89}, and described in details in \cite{sta}.

One needs an ideal $\cq_{\SU}\subset\ker\varepsilon_{\SU}$. The one corresponding to the $4D_+$ calculus is generated by the nine elements $\{c^{2}; \,  c(a^{*}-a); \, q^{2}a^{*2}-(1+q^{2})(aa^{*}-cc^{*})+a^{2}; \, c^{*}(a^{*}-a); \, c^{*2}; \, [q^{2}a+a^{*}-q^{-1}(1+q^{4})]c; \, 
[q^{2}a+a^{*}-q^{-1}(1+q^{4})](a^{*}-a); \, [q^{2}a+a^{*}-q^{-1}(1+q^{4})]c^{*}; \, [q^{2}a+a^{*}-q^{-1}(1+q^{4})][q^{2}a+a^{*}-(1+q^{2})]\}$. One has $\mathrm{Ad}(\cq_{\SU})\subset\cq_{\SU}\otimes\ASU$ and $\dim(\ker\varepsilon_{\SU}/{\cq_{\SU}})=4$. The associated quantum tangent space as in \eqref{qTv} turns out to be a four dimensional $\mathcal{X}_{\cq}\subset\su$.  A choice for a basis  is  given by the elements
\begin{align}
& L_{-}=q^{\frac{1}{2}}FK^{-1},\qquad L_{z}=\frac{K^{-2}-1}{q-q^{-1}},\qquad L_{+}=q^{-\frac{1}{2}}EK^{-1}; \nn \\
& L_{0}=\frac{q(K^{2}-1)+q^{-1}(K^{-2}-1)}{(q-q^{-1})^{2}}\,+FE=
\frac{q(K^{-2}-1)+q^{-1}(K^{2}-1)}{(q-q^{-1})^{2}}\,+EF,
\label{Lq}
\end{align}
from the last commutation rule in \eqref{relsu}. The vector $L_{0}$ belongs to the centre of $\su$: it differs from the quantum Casimir \eqref{cas} by a constant term,
\beq\label{casbis}
C_{q}=L_{0}+\left(\frac{q^{\frac{1}{2}}-q^{-\frac{1}{2}}}{q-q^{-1}}\right)^{2}-\tfrac{1}{4}=L_{0}+[\tfrac{1}{2}]^{2}-\tfrac{1}{4}.
\eeq
The coproducts of the basis \eqref{cpuh} give  $\Delta L_{b}=1\otimes L_{b}+\sum_a L_{a}\otimes f_{ab}$: once chosen the ordering $(-,z,+,0)$, such a tensor product can be represented as a row by column matrix product where
\beq
f_{ab}=\left(\begin{array}{cccc} 1 & 0 & 0 & q^{-\frac{1}{2}}KE \\ 
(q-q^{-1})\qum FK^{-1} & K^{-2} & (q-q^{-1})\qun EK^{-1} & (q-q^{-1})[FE+q^{-1}\,\frac{K^{-2}-K^{2}}{(q-q^{-1})^{2}}] \\ 0 & 0 & 1 & \qun FK \\ 0 & 0 & 0 & K^{2} 
\end{array}
\right).
\label{fab}
\eeq
The differential  $\dd:\ASU\mapsto\Omega^{1}(\SU)$ in  \eqref{ded} is written for any $x\in \ASU$ as
\beq
\dd x=\sum_a (L_{a}\lt x)\omega_{a}=\sum_a \omega_{a}(R_{a}\lt x)
\label{d4}
\eeq
on the dual  basis of left invariant forms $\omega_{a}\in\Omega^{1}(\SU)$ with $\Delta_{L}^{(1)}(\omega_{a})=1\otimes\omega_{a}$.  Here $R_{a}:=-S^{-1}(L_{a})$ and explicitly:
\beq
\label{Rder}
R_{\pm}=L_{\pm}K^2, \quad\quad R_{z}=L_{z}K^2, \quad\quad R_{0}=-L_{0}.
\eeq
On the generators of the algebra the differential  acts as:
\begin{align}
&\dd a=(q-q^{-1})^{-1}(q-1)a\omega_{z}-qc^{*}\omega_{+}+\lambda_{1} a\omega_{0}, \nn \\
&\dd a^{*}=c\omega_{-}+(q-q^{-1})^{-1}(q^{-1}-1)a^{*}\omega_{z}+\lambda_{1} a^{*}\omega_{0}, \nn \\
&\dd c=(q-q^{-1})^{-1}(q-1)c\omega_{z}+a^{*}\omega_{+}+\lambda_{1} c\omega_{0}, \nn \\
&\dd c^{*}=-q^{-1}a\omega_{-}+(q-q^{-1})^{-1}(q^{-1}-1)c^{*}\omega_{z}+\lambda_{1} c^{*}\omega_{0}, 
\label{dcf}
\end{align}  
with $\lambda_1= [\half][\frac{3}{2}]$. 
These relations can be inverted, giving
\begin{align}
&\omega_{-}=c^{*}\dd a^{*}-qa^{*}\dd c^{*}, \qquad \omega_{+}=a\dd c-qc\dd a, \nn \\ 
&\omega_{z}=a^{*}\dd a+c^{*}\dd c-(a\dd a^{*}+q^{2}c\dd c^{*}),\nn \\
&\omega_{0}=(1+q)^{-1}\lambda^{-1}_{1}[a^{*}\dd a+c^{*}\dd c+q(a\dd a^{*}+q^{2}c\dd c^{*})].
\label{om4}
\end{align}
It is then  easy  to see that for $q\to 1$ one has $\omega_{0}\to 0$. This differential calculus reduces in the classical limit  to the standard  three-dimensional bicovariant calculus on $\mathrm{SU(2)}$.

 This first order differential $4D_{+}$ calculus is a $*$-calculus: the $*$-structure on $\ASU$ is  extended to an antilinear  $*$-structure on $\Omega^{1}(\SU)$, such that $(\dd x)^*=\dd (x^*)$ for any $x\in \ASU$. For the basis of left invariant 1-forms is just
\beq
\omega_{-}^{*}=-\omega_{+},\qquad\omega_{z}^{*}=-\omega_{z},\qquad\omega_{0}^{*}=-\omega_{0}.
\label{ss}
\eeq  

From \eqref{bi-struct} one works out the bimodule structure of the calculus, obtaining:
\beq
\begin{array}{lll}
\omega_{-}a=a\omega_{-}-qc^{*}\omega_{0},\qquad &\omega_{+}a=a\omega_{+},\qquad & \omega_{0}a=q^{-1}a\omega_{0},  \\
\omega_{-}a^{*}=a^{*}\omega_{-}, \qquad  &\omega_{+}a^{*}=a^{*}\omega_{+}+c\omega_{0}, \qquad &\omega_{0}a^{*}=qa^{*}\omega_{0},  \\
\omega_{-}c=c\omega_{-}+a^{*}\omega_{0}, \qquad &\omega_{+}c=c\omega_{+}, \qquad & \omega_{0}c=q^{-1}c\omega_{0},  \\
\omega_{-}c^{*}=c^{*}\omega_{-}, \qquad &\omega_{+}c^{*}=c^{*}\omega_{+}-q^{-1}a\omega_{0}, \qquad & \omega_{0}c^{*}=qc^{*}\omega_{0};
\end{array}
\label{biuno}
\eeq
as well as:
\begin{align}
&\omega_{z}a=qa\omega_{z}-q(q-q^{-1})c^{*}\omega_{+}+qa\omega_{0}, \nn \\
&\omega_{z}a^{*}=(q-q^{-1})c\omega_{-}+q^{-1}a^{*}\omega_{z}-q^{-1}a^{*}\omega_{0}, \nn \\
& \omega_{z}c=qc\omega_{z}+(q-q^{-1})a^{*}\omega_{+}+qc\omega_{0}, \nn \\
& \omega_{z}c^{*}=-q^{-1}(q-q^{-1})a\omega_{-}+q^{-1}c^{*}\omega_{z}-q^{-1}c^{*}\omega_{0}.
\label{bidue}
\end{align}

The $\ASU$-bicovariant bimodule $\Omega^{2}(\SU)$ of exterior 2-forms is defined by the projection given in  \eqref{wedk},  with $\mathcal{S}_{\cq}^{(2)}=\ker\,\mathfrak{A}^{(2)}=\ker\,(1-\sigma)\subset\Omega^{1}(\SU)^{\otimes2}$. This necessitates computing the braiding as in \eqref{sigco}, a preliminary step being the computation as in \eqref{ri-co-form} of the right coaction on the left invariant basis forms, $\Delta_{R}^{(1)}(\omega_{a})=\sum_{b}\omega_{b}\otimes J_{ba}$. For the calculus at hand: 
\beq
J_{ba}=\left(\begin{array}{cccc} 
a^{*2} & (1+q^{2})a^{*}c & -qc^{2} & (1-q^{2})a^{*}c \\ 
-qa^{*}c^{*} & aa^{*}-cc^{*} & -ac & (q^{2}-1)cc^{*} \\ 
-qc^{*2} & (q+q^{-1})ac^{*} & a^{2} & (q^{-1}-q)ac^{*} \\ 
0 & 0 & 0 & 1 
\end{array}
\right).
\label{Jba}
\eeq
The braiding map $\sigma:\Omega^{1}(\SU)^{\otimes2}\to\Omega^{1}(\SU)^{\otimes2}$ is then worked out \cite{cla10} to be:
\begin{align}
&\sigma(\omega_{-}\otimes\omega_{-})=\omega_{-}\otimes\omega_{-}, \qquad \sigma(\omega_{+}\otimes\omega_{+})=\omega_{+}\otimes\omega_{+}, \qquad \sigma(\omega_{0}\otimes\omega_{0})=\omega_{0}\otimes\omega_{0}, \nn \\
&\sigma(\omega_{z}\otimes\omega_{z})=\omega_{z}\otimes\omega_{z}+(q^{2}-q^{-2})(\omega_{z}\otimes\omega_{0}+\omega_{-}\otimes\omega_{+}-\omega_{+}\otimes\omega_{-}), \nn \\
&\sigma(\omega_{-}\otimes\omega_{+})=\omega_{+}\otimes\omega_{-}-\omega_{z}\otimes\omega_{0}, \nn \\ 
&\sigma(\omega_{+}\otimes\omega_{-})=\omega_{-}\otimes\omega_{+}+\omega_{z}\otimes\omega_{0}, \nn \\
&\sigma(\omega_{-}\otimes\omega_{z})=\omega_{z}\otimes\omega_{-}+(1+q^{2})\omega_{-}\otimes\omega_{0}, \nn \\
&\sigma(\omega_{z}\otimes\omega_{-})=(1-q^{-2})\omega_{z}\otimes\omega_{-}+q^{-2}\omega_{-}\otimes\omega_{z}-(1+q^{-2})\omega_{-}\otimes\omega_{0},\nn \\
&\sigma(\omega_{-}\otimes\omega_{0})=\omega_{0}\otimes\omega_{-}+(1-q^{2})\omega_{-}\otimes\omega_{0}, \nn \\
&\sigma(\omega_{0}\otimes\omega_{-})=q^{2}\omega_{-}\otimes\omega_{0}, \nn \\
&\sigma(\omega_{z}\otimes\omega_{+})=q^{2}\omega_{+}\otimes\omega_{z}+(1-q^{2})\omega_{z}\otimes\omega_{+}+(1+q^{2})\omega_{+}\otimes\omega_{0}, \nn \\
&\sigma(\omega_{+}\otimes\omega_{z})=\omega_{z}\otimes\omega_{+}-(1+q^{-2})\omega_{+}\otimes\omega_{0}, \nn \\
&\sigma(\omega_{z}\otimes\omega_{0})=\omega_{0}\otimes\omega_{z}+(q-q^{-1})^{2}(\omega_{+}\otimes\omega_{-}-\omega_{-}\otimes\omega_{+})-(q-q^{-1})^{2}\omega_{z}\otimes\omega_{0}, \nn \\
& \sigma(\omega_{0}\otimes\omega_{z})=\omega_{z}\otimes\omega_{0}, \nn \\
&\sigma(\omega_{+}\otimes\omega_{0})=\omega_{0}\otimes\omega_{+}+(1-q^{-2})\omega_{+}\otimes\omega_{0},\nn \\
&\sigma(\omega_{0}\otimes\omega_{+})=q^{-2}\omega_{+}\otimes\omega_{0}.
\label{sig}
\end{align}
Using the general construction of App.\ref{ass:a1}, the $q$-wedge product on 1-forms is defined as $\theta\wedge\theta^{\prime}\,=\,\mathfrak{A}^{(2)}(\theta\otimes\theta^{\prime})\,
=\,(1-\sigma) (\theta\otimes\theta^{\prime})\,\subset\,\mathrm{Range}\,\mathfrak{A}^{(2)}$. On generators: 
\begin{align}
&\omega_{-}\wedge\omega_{-}=\omega_{+}\wedge\omega_{+}=\omega_{0}\wedge\omega_{0}=0, \nn \\
&\omega_{z}\wedge\omega_{z}-(q^{2}-q^{-2})\omega_{+}\wedge\omega_{-}=0, \nn \\
& \omega_{z}\wedge\omega_{\pm}+q^{\pm 2}\omega_{
\pm}\wedge\omega_{z}=0, \nn \\
&\omega_{\pm}\wedge\omega_{0}+\omega_{0}\wedge\omega_{\pm}=0,\nn \\
&\omega_{+}\wedge\omega_{-}+\omega_{-}\wedge\omega_{+}=0, \nn \\
&\omega_{z}\wedge\omega_{0}+\omega_{0}\wedge\omega_{z}-(q-q^{-1})^{2}\omega_{-}\wedge\omega_{+}=0 . 
\label{2fw}
\end{align}
These relations show that $\dim\Omega^{2}(\SU)=6$. The exterior derivative on basis 1-forms results into:
\begin{align}
&\dd\omega_{\pm}= \mp\, q^{\pm 1}\omega_{-}\wedge\omega_{z},  \nn \\
&\dd\omega_{z}=(q+q^{-1})\omega_{+}\wedge\omega_{-}, \nn \\
&\dd\omega_{0}=(q-q^{-1})\,\omega_{-}\wedge\omega_{+}.
\label{d2f}
\end{align}
The  antisymmetriser operator $\mathfrak{A}^{(2)}:\Omega^{2}(\SU)\to\Omega^{2}(\SU)$ 
has a natural spectral decomposition.  This is what we need later on to introduce Hodge operators. A more general analysis of the spectral properties of the antisymmetriser operators associated to a class of bicovariant differential calculi over $\mathrm{SL}_{q}(N)$ (for $N\geq2$) is in \cite{sch99}.  
On the basis 
\beq
\label{p12}
\begin{array}{lll}
\varphi_{0}=\omega_{-}\wedge\omega_{0},& & \varphi_{z}=\omega_{-}\wedge\omega_{0}+(1-q^{-2})\omega_{-}\wedge\omega_{z} \\
\psi_{0}=\omega_{+}\wedge\omega_{0}, & &\psi_{z}=\omega_{+}\wedge\omega_{0}-(1-q^2)\omega_{+}\wedge\omega_{z} \\
\psi_{\pm}=\omega_{0}\wedge\omega_{z}+(1-q^{\pm 2})\omega_{-}\wedge\omega_{+}, & & 
\end{array}
\eeq
which is such that $\varphi_{0}^*=\psi_{0}$, $\varphi_{z}^*=\psi_{z}$ and $\psi_{-}^*=\psi_{+}$, it holds that
\beq
\label{p14}
\begin{array}{lll}
\mathfrak{A}^{(2)}(\varphi_{0})=(1+q^2)\varphi_{0}, &
\mathfrak{A}^{(2)}(\psi_{z})=(1+q^2)\psi_{z}, &
\mathfrak{A}^{(2)}(\psi_{+})=(1+q^2)\psi_{+} \\
\mathfrak{A}^{(2)}(\varphi_{z})=(1+q^{-2})\varphi_{z}, &
\mathfrak{A}^{(2)}(\psi_{0})=(1+q^{-2})\psi_{0}, &
\mathfrak{A}^{(2)}(\psi_{-})=(1+q^{-2})\psi_{-}.
\end{array}
\eeq
For later use we shall adopt the labelling $\xi_{(\pm)}\in\,\mathcal{E}_{(\pm)}$ with 
\beq\label{p12bis}
\mathcal{E}_{(+)}=\{\varphi_{0}, \psi_{z}, \psi_{+}\}, \quad 
\mathrm{and} \qquad \mathcal{E}_{(-)}=\{\varphi_{z},\psi_{0}, \psi_{-}\} .
\eeq

By proceeding further, the $\ASU$-bimodule $\Omega^{3}(\SU)$ is found to be 4-dimensional with left invariant basis elements: 
\beq
\label{p13}
\begin{array}{lll}
\chi_{-}=\omega_{+}\wedge\omega_{0}\wedge\omega_{z},&  &
\chi_{+}=\omega_{-}\wedge\omega_{0}\wedge\omega_{z} \\
\chi_{0}=\omega_{-}\wedge\omega_{+}\wedge\omega_{z}, & & 
\chi_{z}=\omega_{-}\wedge\omega_{+}\wedge\omega_{0},
\end{array}
\eeq
with $\chi_{-}^*=-q^{-2}\chi_{+}$, $\chi_{0}^*=\chi_{0}$ and $\chi_{z}^*=\chi_{z}$.
These exterior forms are closed, 
\beq
\dd\chi_{a}=0, \label{d3f} 
\eeq
and in addition satisfy
\beq
\mathfrak{A}^{(3)}(\chi_{a})=2(1+q^2+q^{-2})\chi_{a} 
\label{p15-}
\eeq
for $a=-,+,z,0$, thus providing the spectral decomposition for the antisymmetriser operator $\mathfrak{A}^{(3)}:\Omega^3(\SU)\to\Omega^3(\SU)$.

The $\ASU$-bimodule $\Omega^{4}(\SU)$ of top forms ($\Omega^{k}(\SU)=\emptyset$ for $k>4$) is 1 dimensional. Its left invariant basis element $\mu=\omega_{-}\wedge\omega_{+}\wedge\omega_{z}\wedge\omega_{0}$ is central, i.e. $x\,\mu=\mu\,x$ for any $x\in\,\ASU$ and its eigenvalue for the action of the antisymmetriser is
\beq
\label{p15}
\mathfrak{A}^{(4)}(\mu)=2(q^4+2q^2+6+2q^{-2}+q^{-4})\mu.
\eeq

\subsection{The exterior algebra over the quantum sphere $\sq$}\label{se:cals2}
The restriction of the first order $4D_+$ bicovariant calculus endows the sphere $\sq$ with a first order left covariant 3-dimensional calculus \cite{ap94,poddc}. 
The exterior algebra $\Omega(\sq)$ can be characterised in terms of some of the  bimodules $\cl_{n}$ introduced in \S\ref{s:qsb}.
Given $f\in \Asq\simeq\cl_{0}$,  the exterior  derivative $\dd:\Asq\mapsto\Omega^{1}(\sq)$ from \eqref{d4} reduces to:  
\beq
\dd f=(L_{-}\lt f)\omega_{-}+(L_{+}\lt f)\omega_{+}+(L_{0}\lt f)\omega_{0}.
\label{d3d}
\eeq
Notice that the basis 1-forms $\{\omega_{a}, a=-,+,0\}$ are graded commutative (cf. \eqref{2fw}). Furthermore, relation \eqref{rellb} shows that $(L_{\pm}\lt f)\in \cl_{\pm2}$ and that $(L_{0}\lt f)\in \cl_{0}$, while the $\ASU$-bimodule structure of $\Omega^{1}(\SU)$ described   by the coproduct \eqref{fab} of the quantum derivations $L_{a}$ gives:
 \beq
\begin{array}{ll}
\phi\,\omega_{-}=\omega_{-}\,\phi-q^{-1}\omega_{0}(L_{+}\lt \phi), & \qquad\qquad \omega_{-}\phi=\phi\,\omega_{-}+q(L_{+}K^{2}\lt \phi)\omega_{0}, \\
\phi^{\prime}\omega_{+}=\omega_{+}\phi^{\prime}-q\omega_{0}(L_{-}\lt\phi^{\prime}), & \qquad\qquad \omega_{+}\phi^{\prime}=\phi^{\prime}\omega_{+}+q^{-1}(L_{-}K^{2}\lt\phi^{\prime})\omega_{0}, \\ 
\phi^{\prime\prime}\omega_{0}=\omega_{0}(K^{-2}\lt\phi^{\prime\prime}), & \qquad\qquad \omega_{0}\phi^{\prime\prime}=(K^{2}\lt\phi^{\prime\prime})\omega_{0}.
\end{array}
\label{bis2}
\eeq
These identities are valid for any $\phi,\phi^{\prime},\phi^{\prime\prime}\in \ASU$. 
They allow to prove by explicit calculations the following identities:
\begin{align}
&\phi\in \cl_{-2}:\qquad\dd(\phi\,\omega_{-})=(L_{+}\lt \phi)\omega_{+}\wedge\omega_{-}+(L_{0}\lt \phi)\omega_{0}\wedge\omega_{-}, \nn \\
&\phi^{\prime}\in \cl_{2}:\qquad \dd(\phi^{\prime}\omega_{+})=(L_{-}\lt\phi^{\prime})\omega_{-}\wedge\omega_{+}+(L_{0}\lt\phi^{\prime})\omega_{0}\wedge\omega_{+}, \nn \\
&\phi^{\prime\prime}\in\cl_{0}:\qquad\dd(\phi^{\prime\prime}\omega_{0})=(L_{-}\lt\phi^{\prime\prime})\omega_{-}\wedge\omega_{0}+(L_{+}\lt\phi^{\prime\prime})\omega_{+}\wedge\omega_{0}+\phi^{\prime\prime}\dd\omega_{0},
\label{deom}
\end{align}
and 
\begin{align}
&\phi\in \cl_{-2}:\qquad\dd(\phi\,\omega_{-}\wedge\omega_{0})=(L_{+}\lt \phi)\omega_{+}\wedge\omega_{-}\wedge\omega_{0}, \nn \\
&\phi^{\prime}\in \cl_{2}:\qquad \dd(\phi^{\prime}\omega_{0}\wedge\omega_{+})=(L_{-}\lt\phi^{\prime})\omega_{-}\wedge\omega_{0}\wedge\omega_{+} , \nn \\
&\phi^{\prime\prime}\in\cl_{0}:\qquad\dd(\phi^{\prime\prime}\omega_{-}\wedge\omega_{+})=(L_{0}\lt\phi^{\prime\prime})\omega_{0}\wedge\omega_{-}\wedge\omega_{+}.  
\label{deom2}
\end{align}
Together with the anti-symmetry properties \eqref{2fw} of the wedge product in $\Omega(\SU)$, these identities suggest that the following proposition holds.
\begin{prop}
\label{lea}
The exterior algebra $\Omega(\sq)$ obtained as a restriction of $\Omega(\SU)$ associated to $4D_{+}$ calculus on $\SU$ can be written in terms of $\Asq$-bimodule isomorphisms:
\begin{align}
\Omega^{1}(\sq)&\simeq\cl_{-2}\,\omega_{-}\oplus\cl_{2}\,\omega_{+}\oplus\cl_{0}\,\omega_{0}
\nn \\
\Omega^{2}(\sq)&\simeq\cl_{-2}\,(\omega_{-}\wedge\omega_{0})\oplus\cl_{0}\,(\omega_{-}\wedge\omega_{+})\oplus\cl_{2}\,(\omega_{0}\wedge\omega_{+}) 
\nn\\
\Omega^{3}(\sq)&\simeq\cl_{0}\,\omega_{-}\wedge\omega_{+}\wedge\omega_{0}
\label{ntd}
\end{align}
\begin{proof}
The analysis above proves only the inclusion $\Omega^{1}(\sq)\subset\cl_{-2}\,\omega_{-}\oplus\cl_{2}\,\omega_{+}\oplus\cl_{0}\,\omega_{0}$ and the analogue ones for higher order forms. The proof of the inverse inclusion will be given at the end of \S\ref{s:ccd}, out of the compatibility of the calculi on the principal Hopf bundle.
\end{proof}
\end{prop}

The basis element $\omega_{-}\wedge\omega_{+}\wedge\omega_{0}$ commutes with all elements in $\cl_{0}\simeq\Asq$.
Such a calculus is 3 dimensional, since from \eqref{d3f} one has $\dd(\phi^{\prime\prime}\omega_{-}\wedge\omega_{+}\wedge\omega_{0})=0$, for any $\phi^{\prime\prime}\in \Asq$, and from \eqref{2fw} one has that $\Omega^{1}(\sq)\wedge(\omega_{-}\wedge\omega_{+}\wedge\omega_{0})=0$.

From \eqref{d4} and \eqref{Rder} the differential can also be written as
\beq
\dd f=\omega_{-}(R_{-}\lt f)+\omega_{+}(R_{+}\lt f)+\omega_{0}(R_{0}\lt f), 
\label{d3di}
\eeq
and it is easy to check the following relations, analogues of the previous \eqref{deom}, \eqref{deom2}:
\begin{align}
&\phi\in \cl_{-2}:\qquad\dd(\omega_{-}\,\phi)=-\omega_{-}\wedge\omega_{+}\,(R_{+}\lt \phi)-\omega_{-}\wedge\omega_{0}\,(R_{0}\lt \phi), \nn \\
&\phi^{\prime}\in \cl_{2}:\qquad \dd(\omega_{+}\,\phi^{\prime})=
-\omega_{+}\wedge\omega_{-}\,(R_{-}\lt \phi^{\prime})-\omega_{+}\wedge\omega_{0}\,(R_{0}\lt \phi^{\prime}), \nn \\
&\phi^{\prime\prime}\in \cl_{0}:\qquad \dd(\omega_{0}\,\phi^{\prime\prime})=
\dd\omega_{0}\wedge\phi^{\prime\prime}-\omega_{0}\wedge\omega_{-}\,(R_{-}\lt \phi^{\prime\prime})-\omega_{0}\wedge\omega_{+}\,(R_{+}\lt \phi^{\prime\prime}); 
\label{deomi}
\end{align}
and 
\begin{align}
&\phi\in \cl_{-2}:\qquad\dd(\omega_{-}\wedge\omega_{0}\,\phi)=\omega_{-}\wedge\omega_{0}\wedge\omega_{+}\,(R_{+}\lt\phi), \nn \\
&\phi^{\prime}\in \cl_{2}:\qquad \dd(\omega_{0}\wedge\omega_{+}\,\phi^{\prime})=\omega_{0}\wedge\omega_{+}\wedge\omega_{-}\,(R_{-}\lt\phi^{\prime}) , \nn \\
&\phi^{\prime\prime}\in\cl_{0}:\qquad\dd(\omega_{-}\wedge\omega_{+}\,\phi^{\prime\prime})=\omega_{-}\wedge\omega_{+}\wedge\omega_{0}\,(R_{0}\lt\phi^{\prime\prime}).  
\label{deom2i}
\end{align}

\section{Hodge operators on $\Omega(\SU)$}
\label{s:Hop}

As described in \S\ref{se:4dc}, it holds  for the bicovariant forms of the $4D_{+}$ first order bicovariant calculus  that the spaces $\Omega^{k}(\SU)$ of forms are free $\ASU$-bimodules with  $\dim\,\Omega^{k}(\SU)=\dim\,\Omega^{4-k}(\SU)$, and $\dim\Omega^4(\SU)=1$. Our  strategy to introduce Hodge operators on 
$\Omega(\SU)$ in \S\ref{sub:H-1} uses first  suitable contraction maps in order to define Hodge operators on the vector spaces $\Omega^{k}_{inv}(\SU)$ of left invariant $k$-forms;  we extend them next to the whole $\Omega^k(\SU)$ by requiring (one side) linearity over $\ASU$.  This follows an alternative although equivalent  
approach to Hodge operators on classical group manifold that we describe first in \S\ref{sse:Hcla}.
A somewhat complementary approach to the one of \S\ref{sub:H-1}, more suitable when restricting to the sphere $\sq$, is then given in \S\ref{sub:H}. 

\subsection{Hodge operators on classical group manifolds}
\label{sse:Hcla}




Let $G$ be an $N$-dimensional  compact connected  Lie group given as a real form of a complex connected 
Lie group. The algebra $\ca(G)=Fun(G)$ of complex valued coordinate functions on $G$ is a  $*$-algebra, whose $*$-structure can be  extended to the whole tensor algebra. A metric on the group $G$ is a non degenerate tensor $g:\mathfrak{X}(G)\otimes\mathfrak{X}(G)\to\ca(G)$ which is symmetric -- i.e. $g(X,Y)=g(Y,X)$, with $X,Y\in\,\mathfrak{X}(G)$ -- and real -- i.e. $g^*(X,Y)=g(Y^*,X^*)$ --. Any metric has a normal form: there exists a basis $\{\theta^{a}, a=1, \dots N\}$ of the $\ca(G)$-bimodule $\Omega^{1}(G)$ of 1-forms which is real, $\theta^{a*}=\theta^a$, such that
\beq
\label{gmf}
g=\sum_{a,b=1}^{N}\eta_{ab}\,\theta^{a}\otimes\theta^b
\eeq
with $\eta_{ab}=\pm1\cdot\delta_{ab}$. 
Given the volume $N$-form $\mu=\mu^*:=\theta^{1}\wedge\ldots\wedge\theta^{N}$, the corresponding Hodge operator $\star:\Omega^{k}(G)\to\Omega^{N-k}(G)$ is the $\ca(G)$-linear operator whose action on the above basis is
\begin{align}
&\star(1)=\mu, \nn \\
\label{Hdui}
&\star(\theta^{a_{1}}\wedge\ldots\wedge\theta^{a_{k}})=\frac{1}{(N-k)!} \sum_{b_j} \epsilon^{a_{1}\ldots a_{k}}_{\qquad b_{1}\ldots b_{N-k}}\theta^{b_{1}}\wedge\ldots\wedge \theta^{b_{N-k}},
\end{align}
with $\epsilon^{a_{1}\ldots a_{k}}_{\qquad b_{1}\ldots b_{N-k}}:=\sum_{s_1\ldots s_k} \eta^{a_{1}s_{1}}\ldots \eta^{a_{k}s_{k}}\epsilon_{s_{1}\ldots s_{k} b_{1}\ldots b_{N-k}}$ from the Levi-Civita tensor and the usual expression for the inverse metric tensor 
$g^{-1}=\sum_{a,b=1}^{N} \eta^{ab} L_{a}\otimes L_{b}$
with  $\sum_{b} \eta^{ab}\eta_{bc}=\delta^{a}_{c}$ on the dual vector field basis such that $\theta^{b}(L_{a})=\delta^{b}_{a}$. The Hodge operator \eqref{Hdui} satisfies the identity:
\beq
\label{quhs}
\star^{2}(\xi)=sgn(g)(-1)^{k(N-k)}\xi
\eeq
on any $\xi\in\,\Omega^k(G)$. Here $sgn(g)=\det(\eta_{ab})$ is the signature of the metric.


Hodge operators can indeed be equivalently introduced in terms of contraction maps. By this we mean an 
$\ca(G)$-sesquilinear map $\Gamma:\Omega^1(G)\times\Omega^1(G)\to\ca(G)$ such that 
$\Gamma(f\,\phi,\eta)=f^*\Gamma(\phi,\eta)$ while $\Gamma(\phi,\eta\,f )=\Gamma(\phi,\eta) f$ for  
$f\in\,\ca(G)$. Such a map can be uniquely extended to a consistent map $\Gamma:\Omega^k(G)\times\Omega^{k+k^\prime}(G)\to\Omega^{k^{\prime}}(G)$. We postpone showing this to the later \S\ref{sub:H-1}
where we prove a similar statement for the bicovariant calculus on $\SU$. 
Having a contraction map, define the tensor $\tilde{g}:\Omega^1(G)\times\Omega^1(G)\to\ca(G)$:  
\beq
\tilde{g}(\phi, \eta):=\Gamma(\phi^{*},\eta).
\label{p7}
\eeq
Next, with a volume form $\mu$, such that $\mu^*=\mu$, define the operator $L:\Omega^{k}(G)\to\Omega^{N-k}(G)$ as
\beq
\label{p3}
L(\xi):=\frac{1}{k!} \Gamma^*(\xi,\mu) 
\eeq
on $\xi\in\,\Omega^k(G)$, having used the notation $\Gamma^*(\cdot,\cdot)=(\Gamma(\cdot,\cdot))^*$.  A second $\ca(G)$-sesquilinear map $\{~,~\}:\Omega^{k}(G)\times\Omega^{k}(G)\to\ca(G)$ can be implicitly introduced by the relation
\beq
\label{p4}
\{\xi,\xi^{\prime}\}\mu:=\xi^*\wedge L(\xi^{\prime}).
\eeq
For any pair of $k$-forms $\xi,\xi^{\prime}$ it is straightforward  to recover that
\beq
\label{pi5}
\{\xi,\xi^{\prime}\}=\frac{1}{k!} \Gamma^*(\xi^{\prime},\xi) .
\eeq
The operator \eqref{p3} is not in general an Hodge operator: one has for example $L(1)=\mu$ as well as $L(\mu)= \det (\Gamma^*(\mu,\mu))$ which is not necessarily $\pm 1$. 
To recover the standard formulation for a Hodge operator, one has to impose two constraints:
\begin{enumerate}[(a)]
\item An hermitianity condition. The sesquilinear map $\Gamma$ is said hermitian provided it satisfies:
\beq
\label{p10}
\{\phi,\eta\}=\Gamma(\phi,\eta),
\eeq
for any couple of 1-forms $\phi$ and $\eta$. 
\end{enumerate}
From \eqref{pi5} and \eqref{p4} it holds  that    $\{\phi,\eta\}=\Gamma^*(\eta,\phi)$.Then  
\beq
\label{pi10}
\{\phi,\eta\}=\Gamma(\phi,\eta)\quad\Leftrightarrow\quad\Gamma(\phi,\eta)=\Gamma^{*}(\eta,\phi).
\eeq
If the sesquilinear form $\Gamma$ is hermitian, one can prove that the expression \eqref{pi5} becomes 
\beq
\label{pii10}
\{\xi,\xi^{\prime}\}=\frac{1}{k!}\Gamma(\xi,\xi^{\prime}).
\eeq
\begin{enumerate}[(b)]
\item A reality condition, namely a compatibility of the operator $L$ with the $*$-conjugation:
\beq
\label{p11}
L(\phi^{*})=(L(\phi))^*
\eeq
on 1-forms.
\end{enumerate}
If these two constraints are fullfilled, the tensor $\tilde{g}$ in \eqref{p7} is symmetric and real: it is (the inverse of) a metric tensor on the group manifold $G$. The  operator $L$ turns out to be the standard Hodge operator corresponding to the metric given by $\tilde{g}$, and satisfies the identities:
\begin{align}
\label{p8}
L^2(\xi)=(-1)^{k(N-k)}sgn(\Gamma)\xi , \qquad
\{\xi,\xi^{\prime}\}= sgn(\Gamma)\{L(\xi),L(\xi^{\prime})\} 
\end{align}
with 
$$sgn(\Gamma):=(\det(\Gamma(\phi^{a},\phi^{b}))|\det(\Gamma(\phi^{a},\phi^{b})|^{-1}=sgn(\tilde{g}).
$$
Moreover, the operator $L$ turns out to be real, that is, it commutes with the hermitian conjugation $*$, on the whole exterior algebra $\Omega(G)$. 

The above procedure could be somehow inverted. That is, given an hermitian contraction map $\Gamma$ as in \eqref{p10},  define the  operator $L$ by \eqref{p3}. The corresponding tensor $\tilde{g}$ turns out to be  real, but non necessarily symmetric. Imposing  $L$ to satisfy one of the two conditions in \eqref{p8} -- they are proven to be  equivalent --  makes the tensor $\tilde{g}$ symmetric, that is the inverse of a metric tensor, whose Hodge operator is  $L$. 
 
\subsection{Hodge operators on $\Omega(\SU)$}\label{sub:H-1}

In this section we shall describe how the classical geometry analysis of the previous section can be used 
to introduce an Hodge operator on both the exterior algebras $\Omega(\SU)$ and $\Omega(\sq)$ built out of the $4D$-bicovariant calculus \`a la Woronowicz on $\ASU$. 
A somewhat different formulation of contraction maps was also used 
in \cite{hec99,hec03} for a family of Hodge operators on the exterior algebras of bicovariant differential calculi over quantum groups.

We shall then start with a contraction map $\Gamma:\Omega^1_{inv}(\SU)\times\Omega^1_{inv}(\SU)\to\IC$,
required to satisfy $\Gamma(\lambda\,\omega,\omega^{\prime}) = \lambda^* \Gamma(\omega, \omega^{\prime})$ and 
$\Gamma(\omega,\omega^{\prime}\, \lambda) = \Gamma(\omega, \omega^{\prime}) \lambda$, for $\lambda\in\IC$.  
The natural extension to $\Gamma:\Omega_{inv}^{\otimes k}(\SU)\times\Omega_{inv}^{\otimes k+k^{\prime}}(SU)\to\Omega_{inv}^{\otimes k^{\prime}}(SU)$ given by
\beq
\label{p2i}
\Gamma(\omega_{a_{1}}\otimes\ldots\otimes\omega_{a_k},\omega_{b_1}\otimes\ldots\otimes\omega_{b_{k+k^{\prime}}}):=\left(\Pi_{j=1}^{k}\,\Gamma(\omega_{a_{j}},\omega_{ b_{j}})\right)\,\omega_{b_{k+1}}\otimes\cdots\otimes\omega_{b_{k+k^{\prime}}} ,
\eeq
with the assumption that $\Gamma(1,\omega)=\omega$ for any $\omega\in\,\Omega(\SU)$, can be used to define a consistent  contraction map 
$\Gamma:\Omega^{ k}(\SU)\times\Omega^{ k+k^{\prime}}(\SU)\to\Omega^{ k^{\prime}}(\SU)$, via 
\begin{multline}
\label{c0}
\Gamma(\omega_{a_{1}}\wedge\ldots\wedge\omega_{a_k},\omega_{b_1}\wedge\ldots\wedge\omega_{b_{k+k^{\prime}}}) \\
:=\Gamma(\mathfrak{A}^{(k)}(\omega_{a_{1}}\otimes\ldots\otimes\omega_{a_k}), \mathfrak{A}^{(k+k^{\prime})}(\omega_{b_1}\otimes\ldots\otimes\omega_{b_{k+k^{\prime}}})) .
\end{multline}
This comes from the $k$-th order anti-symmetriser $\mathfrak{A}^{(k)}$, constructed from the braiding of the calculus, and used to define the exterior product of forms,
\beq
\label{1c}
\omega_{a_{1}}\wedge\ldots\wedge\omega_{a_{k}}:=\mathfrak{A}^{(k)}(\omega_{a_{1}}\otimes\cdots\otimes\omega_{a_{k}});
\eeq
the key identity for the consistency of \eqref{c0} is 
\beq
\label{2c}
	\mathfrak{A}^{(k+k^{\prime})}(\omega_{a_{1}}\otimes\cdots\otimes\omega_{a_{k+k^{\prime}}})=(\mathfrak{A}^{(k)}\otimes\mathfrak{A}^{(k^{\prime})})(\sum_{\sigma_j\in\mathit{S}(k,k^{\prime})}(-1)^{\pi_{\sigma_{j}}}\sigma_j(\omega_{a_{1}}\otimes\cdots\otimes\omega_{a_{k+k^{\prime}}})),
\eeq
where $\mathit{S}(k,k^{\prime})$ is the collection of the $(k,k^{\prime})$-shuffles, permutations $\sigma_{j}$ of $\{1,\ldots,k+k^{\prime}\}$ such that $\sigma_{j}(1)<\cdots<\sigma_{j}(k)$ and $\sigma_{j}(k+1)<\ldots<\sigma_{j}(k+k^{\prime})$, and $\pi_{\sigma_{j}}$ is the parity of $\sigma_{j}$.  The identity \eqref{2c} is valid on the whole exterior algebra over any bicovariant calculus \`a la Woronowicz on a quantum group. It allows  to show that any $(k+k^{\prime})$-form can be written as a linear combination of  tensor products of  $k$-forms times  $k^{\prime}$-forms.

To proceed further, we use a slightly more general volume form by taking $\mu=\mu^*= \ii \, m\,\omega_{-}\wedge\omega_{+}\wedge\omega_{0}\wedge\omega_{z}$, with $m\in\,\IR$. Then we define an operator 
$$
\star:\Omega^{k}_{inv}(\SU)\to\Omega^{4-k}_{inv}(\SU),
$$
in degree zero and one by  
\beq\label{p17}
\star(1):= \Gamma^*(1,\mu) =\mu  \qquad \mathrm{and} \qquad 
\star(\omega_{a}):= \Gamma^*(\omega_{a},\mu) . 
\eeq
For $\Omega^{k}_{inv}(\SU)$ with $k\geq2$ we use the diagonal bases of the antisymmetriser, that is 
\beq
\mathfrak{A}^{(k)}(\xi)=\lambda_{\xi}\xi , 
\label{eixi}
\eeq
with coefficients in \eqref{p14}, \eqref{p15-} and \eqref{p15} respectively. On these basis we define 
\beq
\star(\xi):=\frac{1}{\lambda_{\xi}}\, \Gamma^*(\xi,\mu) .
\label{p17bis}
\eeq
Here and in the following we denote $(\Gamma(~,~))^*=\Gamma^*(~,~)$.
The definition \eqref{p17bis} is a natural generalisation of the classical \eqref{p3}: the classical factor $k!$ -- the spectrum of the antisymmetriser operator on $k$-forms in the classical case, where the braiding is the flip operator --  is replaced by the spectrum of the quantum antisymmetriser. Also, the presence of the *-conjugate comes from consistency and in order to have non trivial solutions. 

Before we proceed, it is useful to re-express the volume forms in terms of the diagonal bases of the anti-symmetriser  operators. Some little algebra shows that 
\begin{align}\label{p16}
&\mu=\ii m \{-\omega_{-}\otimes\chi_{+}^*+\omega_{+}\otimes\chi_{-}^*+\omega_{0}\otimes\chi_{0}^*-\omega_{z}\otimes\chi_{z}^*\}  \nn \\
&\mu=\ii m \{-\chi_{z}\otimes\omega_{z}^*+\chi_{-}\otimes\omega_{+}^* -\chi_{+}\otimes\omega_{-}^* +\chi_{0}\otimes\omega_{0}^*\}, 
\end{align} 
and 
\begin{multline}\label{p16bis}
\mu=\frac{\ii m}{q^2-1}\Big\{\frac{1}{1+q^2}(q^4\psi_{-}\otimes\psi_{+}^*-\psi_{+}\otimes\psi_{-}^*)  \\ 
  +(q^4\varphi_{z}\otimes\varphi_{0}^*-\varphi_{0}\otimes\varphi_{z}^*+q^2\psi_{0}\otimes\psi_{z}^*-q^{-2}\psi_{z}\otimes\psi_{0}^*) \Big\} . 
\end{multline}
A little more algebra shows in turn that on 1-forms
\beq\label{p19-1}
\star(\omega_{a})=\ii m  \big\{ \Gamma^*(\omega_{a},\omega_{-}) \chi_{+} 
- \Gamma^*(\omega_{a},\omega_{+}) \chi_{-} 
- \Gamma^*(\omega_{a},\omega_{0}) \chi_{0} 
+ \Gamma^*(\omega_{a},\omega_{z}) \chi_{z} 
\big\} ; 
\eeq
and that using the bases \eqref{p12bis}, on 2-forms
\begin{multline}\label{p19-2}
 \star(\xi_{(+)})=\frac{\ii m}{1-q^4} \Big\{
\frac{1}{1+q^2} \left(q^4\Gamma^*(\xi_{(+)},\psi_{-})\psi_{+}-\Gamma^*(\xi_{(+)},\psi_{+})\psi_{-} \right)+(q^4\Gamma^*(\xi_{(+)},\varphi_{z})\varphi_{0}  \\
-\Gamma^*(\xi_{(+)},\varphi_{0})\varphi_{z}+q^2\Gamma^*(\xi_{(+)},\psi_{0})\psi_{z}-q^{-2}\Gamma^*(\xi_{(+)},\psi_{z})\psi_{0}) \Big\} 
\end{multline}
\begin{multline*}
 \star(\xi_{(-)})=\frac{\ii m}{q^{-2}-q^{2}} \Big\{
\frac{1}{1+q^2}  \left( q^4\Gamma^*(\xi_{(-)},\psi_{-})\psi_{+}-\Gamma^*(\xi_{(-)},\psi_{+})\psi_{-} \right)+(q^4\Gamma^*(\xi_{(-)},\varphi_{z})\varphi_{0}    \\
 -\Gamma^*(\xi_{(-)},\varphi_{0})\varphi_{z}+q^2\Gamma^*(\xi_{(-)},\psi_{0})\psi_{z}-q^{-2}\Gamma^*(\xi_{(-)},\psi_{z})\psi_{0})
\Big\}. 
\end{multline*}
As for 3-forms one finds
\begin{multline}\label{p19-3}
 \star(\chi_{a})=-\frac{\ii m}{2(1+q^2+q^{-2})} \Big\{-\Gamma^*(\chi_{a},\chi_{+})\omega_{-} \\ +\Gamma^*(\chi_{a},\chi_{-})\omega_{+}+\Gamma^*(\chi_{a},\chi_{0})\omega_{0}-\Gamma^*(\chi_{a},\chi_{z})\omega_{z} \Big\},   \end{multline}
and finally for the top form
\beq\label{p19-4}
 \star(\mu)=\frac{1}{2(q^4+2q^2+6+2q^{-2}+q^{-4})} \Gamma^*(\mu,\mu).
\eeq

As in \eqref{p4} we define the sesquilinear map $\{~,~\}:\Omega^{k}_{inv}(\SU)\times\Omega^{k}_{inv}(\SU)\to\IC$ by
\beq
\{\xi,\xi^{\prime}\}\mu:=\xi^*\wedge\star(\xi^{\prime}).
\label{p4q}
\eeq
Then, mimicking the analogous construction of \S\ref{sse:Hcla}  we impose both an hermitianity and a reality condition on the contraction  map. 
\begin{enumerate}[(a)]
\item A contraction map  is hermitian provided it satisfies:
\beq
\label{pi21}
\{\omega_{a},\omega_{b}\}=\Gamma(\omega_{a},\omega_{b}), \qquad \mathrm{for} \qquad  a,b=-,+,z,0 .
\eeq 
\end{enumerate}
Given contraction maps fullfilling such an hermitianity constraint, from the first line in \eqref{p19-1} one has 
that $\Gamma(\omega_{a},\omega_{b}) = \Gamma^*(\omega_{b},\omega_{a})$. i.e. $\Gamma_{ab}=\Gamma^{*}_{ba}$. With such a condition it is moreover possible to prove, that for with $k=2,3,4$,
\beq
\label{p21}
\{\xi,\xi^{\prime}\}=\frac{\lambda_{\xi^*}}{\lambda_{\xi}\lambda_{\xi^{\prime}}}\,\Gamma(\xi,\xi^{\prime}).
\eeq
on any  $\xi,\xi^{\prime}\in\,\Omega^{k}_{inv}(\SU)$ of a diagonal basis of the antisymmetrizer as in \eqref{eixi}.
The above expression is the counterpart of \eqref{pii10} for a braiding which is not just the flip operator. 
\begin{enumerate}[(b)]
\item
An hermitian contraction map is real provided one has
\beq
\label{pi22}
\lambda_{\xi^*}(\star\xi^*)=(\lambda_{\xi}(\star\xi))^*.
\eeq
\end{enumerate}
again on a diagonal basis of $\mathfrak{A}^{(k)}(\xi)$.
 This expression generalises the classical one \eqref{p11}. Notice that it is set on any  $\Omega^{k}_{inv}(\SU)$, and not only on 1-forms as in the classical case.   

The requirement that the contraction be hermitian and real results in a series of constraints. Firstly, 
the action on $\Omega^1_{inv}(\SU)$ of the corresponding operator $\star$ as defined in \eqref{p17} is worked out to be given by
\beq
 \label{piii22}
 \star\left(\begin{array}{c} \omega_{-} \\ \omega_{+} \\ \omega_{0} \\ \omega_{z} \end{array}\right)=\ii m \left(\begin{array}{cccc}  0 & \alpha & 0 & 0 \\ -q^2\alpha  & 0 &  0 & 0 \\ 0 & 0 &  -\nu & \epsilon \\ 0 & 0 & -\epsilon & \gamma \end{array} \right)\left(\begin{array}{c}\chi_{-} \\ \chi_{+} \\ \chi_{0} \\ \chi_{z} \end{array}\right).
 \eeq
The only non zero terms of the contraction $\Gamma$ are given by
\begin{align}
\Gamma_{--}=q^{-2}\Gamma_{++}=\alpha, \qquad \Gamma_{0z}=\Gamma_{z0}=\epsilon, \
\qquad \Gamma_{00}=\nu, \qquad
\Gamma_{zz}=\gamma, 
\label{piv22}
\end{align}
with parameters that are real and satisfy in addition the conditions:
\begin{align}
&2\nu+(q^2-q^{-2})\epsilon=0, \nn \\
&2(\epsilon^2-\gamma\nu)+(q-q^{-1})^2(2q^2\alpha^2+\epsilon^2)=0.
\label{pv22}
\end{align} 

On $\Omega^{2}_{inv}(\SU)$  the action of such  operator is block off-diagonal,  
\begin{align}
 \label{pvi22}
 &\star\left(\begin{array}{c}\varphi_{0} \\ \psi_{z} \\ \psi_{+}  \end{array}\right)\,= 
 \frac{\ii m}{q^4-1}\left(\begin{array}{ccc} \Gamma(\varphi_{0},\varphi_{0}) & 0 & 0 \\ 
  0 & q^4\Gamma(\varphi_{z},\varphi_{z}) & 0 \\  0 & 0 & \frac{\Gamma(\psi_{+},\psi_{+})}{1+q^{2}} 
 \end{array} 
 \right)\left(\begin{array}{c} \varphi_{z} \\ \psi_{0} \\ \psi_{-} \end{array}\right),\nn
 \\  
 ~\nn \\
&\star\left(\begin{array}{c} \varphi_{z} \\ \psi_{0} \\ \psi_{-} \end{array}\right)\,= 
 \frac{\ii m}{1-q^4}\left(\begin{array}{ccc}  q^6\Gamma(\varphi_{z},\varphi_{z}) & 0 & 0  \\ 
0  & q^2\Gamma(\varphi_{0},\varphi_{0}) & 0  \\ 0 & 0 & \frac{q^4\Gamma(\psi_{+},\psi_{+})}{1+q^{2}}
 \end{array} 
 \right)\left(\begin{array}{c}\varphi_{0} \\ \psi_{z} \\ \psi_{+}  \end{array}\right),
\end{align}

\noindent
while on $\Omega^{3}_{inv}(\SU)$ is
 \begin{multline}
 \label{pvii22}
 \star\left(\begin{array}{c} \chi_{-} \\ \chi_{+} \\ \chi_{0} \\ \chi_{z} \end{array}\right)= \\ =
 \frac{\ii m}{2(1+q^2+q^{-2})}\left(\begin{array}{cccc}0 & -\Gamma(\chi_{-},\chi_{-}) & 0 & 0  \\  q^2\Gamma(\chi_{-},\chi_{-}) & 0 & 0 & 0 \\ 0 & 0 & -\Gamma(\chi_{0},\chi_{0}) & \Gamma(\chi_{z},\chi_{0}) \\ 0 & 0 & -\Gamma(\chi_{0},\chi_{z}) & \Gamma(\chi_{z},\chi_{z}) \end{array}\right)
 \left(\begin{array}{c} \omega_{-} \\ \omega_{+} \\ \omega_{0} \\ \omega_{z} \end{array}\right).
 \end{multline}

\bigskip
It turns out that the square of the operator $\star$ is not necessarily diagonal. 
An explicit computation shows moreover that when $q\neq1$, given the constraints \eqref{pv22} there is no choice for the contraction $\Gamma$, nor for the value of the scale parameter $m\in\,\IR$ in the volume form such that the spectrum of the operator $\star^2$ is constant on any vector space $\Omega^k_{inv}(\SU)$. 
This means that the operator $\star$ does not satisfy the classical expressions in \eqref{p8}. 
We choose a particular value for the parameter $m$ defining
\begin{align}
\det\Gamma:=\frac{1}{\lambda_{\mu}}\Gamma(\omega_{-}\wedge\omega_{+}\wedge\omega_{0}\wedge\omega_{z},\omega_{-}\wedge\omega_{+}\wedge\omega_{0}\wedge\omega_{z}), \qquad  
sgn(\Gamma):=\frac{\det\Gamma}{|\det\Gamma|}
\label{pviii22}
\end{align}
and imposing
\beq
\label{pix22}
\star^2(1)=sgn(\Gamma),
\eeq
which is clearly equivalent to the constraint 
\beq
\label{px22}
m^2=|\det\Gamma|^{-1}.
\eeq
An explicit calculation shows that conditions \eqref{pi21} and \eqref{pi22} fix the quantum determinant \eqref{pviii22}  to be positive, so that we have $sgn(\Gamma)=1$.

We finally extend the operator $\star$ to the whole exterior algebra. This can be defined in two ways, i.e. we  define Hodge operators $\star^{L},\star^{R}:\Omega^{k}(\SU)\to\Omega^{4-k}(\SU)$ by: 
\begin{align}
 \star^{L}(x\,\omega):=x\star(\omega), \qquad\qquad\qquad
 \star^{R}(\omega\,x):=(\star\, \omega) x ,
\label{Hlr}
\end{align}
with $x\in\,\ASU$ and $\omega\in\,\Omega_{inv}(\SU)$. Both operators will find their use later on.

\subsection{Hodge operators on $\Omega(\SU)$ -- a complementary approach} \label{sub:H}

The procedure used in the previous section cannot be extended \emph{ipso facto} to introduce an Hodge operator on the exterior algebra $\Omega(\sq)$: although all $\Omega^{k}(\sq)$ are free left $\Asq$-modules \cite{HK03},  the tensor product  $\Omega^{\otimes2}(\sq)$ has no braiding like the $\sigma$ above.

In order to construct a suitable Hodge operator  on the quantum sphere, we shall  export to this quantum homogeneous space the construction of \cite{kmt},  originally conceived on the exterior algebra over a quantum group.  
The strategy largely coincides with the one described in \cite{ale09} and presents similarities to that used 
in \cite{dal09} where a Hodge operator has been introduced on a quantum projective plane.

We start by briefly recalling the formulation from \cite{kmt}. 
Consider a  $*$-Hopf algebra $\ch$ and the exterior algebra $\Omega(\ch)$ over  an $N$-dimensional left covariant first order calculus 
$(\Omega^{1}(\ch),\dd)$,  with  $\dim \Omega^{N-k}(\ch)=\dim \Omega^{k}(\ch)$ and $\dim \Omega^{N}(\ch)=1$. Suppose in addition that $\ch$ has an Haar state $h:\ch\to\IC$, i.e. a unital functional, which is 
 invariant, i.e.  $(\id\otimes h)\Delta x=(h\otimes \id)\Delta x=h(x)1$ for any $x\in \ch$, and positive,  i.e.  $h(x^{*}x)\geq0$ for all $x\in \ch$.  An Haar state so defined  is unique and automatically faithful: $h(x^{*}x)=0$ implies $x=0$. 
Upon fixing an inner product on a left invariant basis of forms, the state $h$ is then used to endow 
the whole exterior algebra with a left and a right inner product, when requiring left or right invariance,
\begin{align}
&\hs{x\,\omega}{x^{\prime}\,\omega^{\prime}}^{L}:=h(x^{*}x^{\prime})\hs{\omega}{\omega^{\prime}} , 
\nn \\
&\hs{\omega\,x}{\omega^{\prime}\,x^{\prime}}^{R}:=h(x^*x^{\prime})\hs{\omega}{\omega^{\prime}} 
\label{inpo}
\end{align}
for any $x,x^{\prime}\in \ch$ and  $\omega,\omega^{\prime}$ in $\Omega_{inv}(\ch)$.  
The spaces $\Omega^{k}(\ch)$ are taken to be pairwise orthogonal (this is stated by saying that 
the inner product is graded). 

The differential calculus is said to be non-degenerate if, whenever $\eta\in \Omega^{k}(\ch)$ and $\eta^{\prime}\wedge\eta=0$ for any $\eta^{\prime}\in \Omega^{N-k}(\ch)$, then necessarily $\eta=0$. 
Choose in $\Omega^{N}(\ch)$ a left invariant hermitian basis element $\mu=\mu^{*}$, referred to as the volume form of the calculus. For the sake of the present paper, we assume that the differential calculus has a volume form such that $\mu\,x=x\,\mu$ for any $x\in\,\ch$ (this condition is satisfied by the $4D_{+}$ bicovariant calculus on $\SU$ that we are considering). Then one defines an `integral'
\begin{align*}
&\int _{\mu} \, : \, \Omega(\ch)\to\IC, \qquad\qquad 
\int _{\mu} x\, \mu =h(x) , \qquad \mathrm{for} \quad x\in \ch ,
\end{align*}
and $\int_{\mu} \eta=0$ for any $k$-form 
$\eta$ with $k<N$. For a non-degenerate calculus the functional  $\int_{\mu}$  is left-faithful if $\eta \in \Omega^{k}(\ch)$ is such that $\int_{\mu}\eta^{\prime}\wedge\eta=0$ for all $\eta^{\prime}\in \Omega^{N-k}(\ch)$, then $\eta=0$.  
The central result is \cite{kmt}: 
\begin{prop}
\label{Lge}
Consider a left covariant, non-degenerate differential calculus on a $*$-Hopf algebra, whose corresponding exterior algebra is 
such that $\dim\, \Omega^{N-k}(\ch)=\dim\, \Omega^{k}(\ch)$ and $\dim\, \Omega^{N}(\ch)=1$, with a left-invariant volume form $\mu=\mu^*$ satisfying $x\,\mu=\mu\,x$ for any $x\in\,\ch$.  If $\Omega(\ch)$ is endowed with  inner products and  integrals  as before, there exists a unique left $\ch$-linear bijective operator 
$L:\Omega^{k}(\ch)\to\Omega^{N-k}(\ch)$ for $k=0,\ldots,N$  (resp.  a unique right $\ch$-linear bijective operator $R$) such that 
\beq\label{Lop} 
\int_{\mu}\eta^{*}\wedge L(\eta^{\prime})=\hs{\eta}{\eta^{\prime}}^L, \qquad\qquad  
\int_{\mu}\eta^{*}\wedge R(\eta^{\prime})=\hs{\eta}{\eta^{\prime}}^R
\eeq 
for any $\eta,\eta^{\prime}\in \Omega^{k}(\ch)$. 
\end{prop}
We mention that there is no $R$ operator in \cite{kmt}. It is just to prove its right $\ch$-linearity that one needs the condition $x\mu=\mu x$ for the volume form $\mu$ with $x\in\,\ch$.

We are now ready to make contact with the previous \S\ref{sub:H-1}. The $4D_{+}$ differential calculus on 
$\SU$ is easily seen to be non degenerate. On the other hand, the Haar state functional $h$ is given by (cf. \cite{KS97}): 
\beq
h(1)=1;\qquad h((cc^{*})^{k})=(\sum_{j=0}^{k}\,q^{2j})^{-1}=\frac{1}{1+q^{2}+\ldots+q^{2k}},
\label{Has}
\eeq
with $k\in \IN$,  all other generators mapping to zero.

Now, use the sesquilinear map \eqref{p4q}  for an inner product 
$\hs{\omega}{\omega^{\prime}}:=\{\omega,\omega^{\prime}\}$ on generators of $\Omega_{inv}(\SU)$ and extend it to a left invariant and a right invariant ones to the whole of $\Omega_{inv}(\SU)$ as in \eqref{inpo} using the state $h$. The uniqueness of the operators $L$ and $R$ from Proposition~\ref{Lge} then implies that the extended left and right  inner products are related to the left and right Hodge operators \eqref{Hlr}  by 
\beq
\label{Holr}
\int_{\mu}\eta^{*}\wedge (\star^{L}\,\eta^{\prime})=\hs{\eta}{\eta^{\prime}}^L, \qquad\qquad  
\int_{\mu}\eta^{*}\wedge (\star^{R}\,\eta^{\prime})=\hs{\eta}{\eta^{\prime}}^R
\eeq 
for any $\eta,\eta^{\prime}\in \Omega^{k}(\ch)$. 

\section{Hodge operators on $\Omega(\sq)$} \label{s:hl3}
From the previous section, 
the procedure to introduce Hodge operators on the quantum sphere appears outlined. Inner products on $\Omega(\SU)$ naturally induce inner products on $\Omega(\sq)$, and we shall explore the use of relations like the \eqref{Holr} above to define a class of Hodge operators.   


The exterior algebra $\Omega(\sq)$ over the quantum sphere $\sq$ is described in \S\ref{se:cals2}. In particular, we recall  its description  in terms of the $\Asq$-bimodules $\cl_{n}$ given in \eqref{libu}:
\begin{align}
&\Omega^{0}(\sq) \simeq \Asq \simeq \cl_{0} \,  \nn \\
&\Omega^{1}(\sq) \simeq \cl_{-2}\,\omega_{-}\oplus\cl_{2}\,\omega_{+}\oplus\cl_{0}\,\omega_{0}\,\simeq \omega_{-}\,\cl_{-2}\oplus\omega_{+}\,\cl_{2}\oplus \omega_{0}\,\cl_{0} \nn \\
&\Omega^{2}(\sq) \simeq \cl_{-2}\,(\omega_{-}\wedge\omega_{0})\oplus\cl_{0}\,(\omega_{-}\wedge\omega_{+})\oplus\cl_{2}\,(\omega_{0}\wedge\omega_{+})\nn \\
&\qquad\qquad\qquad\qquad \simeq  
(\omega_{-}\wedge\omega_{0})\,\cl_{-2}\oplus(\omega_{-}\wedge\omega_{+})\,\cl_{0}\oplus(\omega_{0}\wedge\omega_{+})\,\cl_{2}  \,  \nn \\
&\Omega^{3}(\sq) \simeq \cl_{0}\,\omega_{-}\wedge\omega_{+}\wedge\omega_{0} \label{isoc2}  \,\simeq\,\omega_{-}\wedge\omega_{+}\wedge\omega_{0}\,\cl_{0} \, .
\end{align} 
In the rest of this section, to be consistent with the notation introduced in \S\ref{se:cals2}, we shall consider  elements $\phi,\psi\in\,\cl_{-2}$, elements  $\phi^{\prime},\psi^{\prime}\in\,\cl_{2}$ and elements  $\phi^{\prime\prime},\psi^{\prime\prime}\in\,\cl_{0}$.
\begin{lemm}
The above left covariant 3D calculus on $\sq$  is non-degenerate.
\label{nddc}
\begin{proof}
Given $\theta\in \Omega^{k}(\sq)$ the condition of non degeneracy, namely  $\theta^{\prime}\wedge\theta=0$ for any $\theta^{\prime}\in \Omega^{3-k}(\sq)$ only if $\theta=0$,  is trivially satisfied for $k=0,3$.
  
From \eqref{isoc2}  take the 1-form 
$\theta=\phi\,\omega_{-}$  and  a 2-form $\theta^{\prime}=\psi\,\omega_{-}\wedge\omega_{0}+\psi^{\prime}\omega_{+}\wedge\omega_{0}+\psi^{\prime\prime}\omega_{-}\wedge\omega_{+}$. Using the commutation properties \eqref{bis2} between 1-forms and elements in $\ASU$,  one has $\theta^{\prime}\wedge\theta=\{\psi^{\prime}(K^{2}\lt\phi)-\psi^{\prime\prime}(\qun KE\lt\phi)\}\omega_{-}\wedge\omega_{+}\wedge\omega_{0}$, so that the equation  $\theta^{\prime}\wedge\theta=0$ for any $\theta^{\prime}\in \Omega^{2}(\sq)$ is equivalent to the condition $\{\psi^{\prime}(K^{2}\lt\phi)-\psi^{\prime\prime}(\qun KE\lt\phi)\}=0$ for any $\theta^{\prime}\in \Omega^{2}(\sq)$; taking $\theta^{\prime}=\psi^{\prime}\omega_{+}\wedge\omega_{0}$, one shows that this condition is satisfied only if $\phi=0$. A similar conclusion is reached with a 1-form $\theta=\phi^{\prime}\omega_{+}$, and with a 1-form $\theta=\phi^{\prime\prime}\omega_{0}$. 

Consider then a 2-form $\theta=\phi\,\omega_{-}\wedge\omega_{0}$, and a 1-form $\theta^{\prime}=\psi\omega_{-}+\psi^{\prime}\omega_{+}+\psi^{\prime\prime}\omega_{0}$. Their product  is $\theta^{\prime}\wedge\theta=(\psi^{\prime}\phi)\omega_{+}\wedge\omega_{-}\wedge\omega_{0}$, so that the condition $\theta^{\prime}\wedge\theta=0$, for all $\theta^{\prime}\in \Omega^{1}(\sq)$ is equivalent to the condition $\psi^{\prime}\phi=0$ for any $\psi^{\prime}$; this condition is obviously satisfied only by $\phi=0$. It is clear that a similar analysis can be performed for  any 2-form $\theta\in \Omega^{2}(\sq)$.
\end{proof}
\end{lemm}

The  Haar state $h$ of $\ASU$ given in \eqref{Has} yields a faithful and invariant state when restricted to 
$\Asq$.  As a volume form we take $\check{\mu}=\check{m}\,\omega_{-}\wedge\omega_{+}\wedge\omega_{0}=\check{\mu}^*$ with $\check{m}\in \IR$. It commutes with every  algebra element, $f\,\check{\mu}=\check{\mu}\,f$ for $f\in\,\Asq$, so the  integral  on the exterior  algebra $\Omega(\sq)$ can be defined by
\begin{align}
&\int_{\check{\mu}}\theta=0,\qquad &\mathrm{on}\quad \theta\in \Omega^{k}(\sq),\,\mathrm{for}\,k=0,1,2 \, , 
\nn \\
&\int_{\check{\mu}}f\,\check{\mu}=h(f) ,\qquad &\mathrm{on}\quad f\,\check{\mu} \in \Omega^{3}(\sq) \, . 
\label{ints2}
\end{align}

\begin{lemm}
\label{lef}
The integral $\int_{\check{\mu}}:\Omega(\sq)\to\IC$ defined by \eqref{ints2} is left-faithful.
\begin{proof}
The  proof of the left-faithfulness of the integral can be easily established from a direct analysis,  using the faithfulness of the  Haar state $h$. \end{proof}
\end{lemm}

The restriction to $\Omega(\sq)$ of the left and right $\ASU$-linear graded inner products on $\Omega(\SU)$ in \eqref{Holr}  gives  left  and right 
$\Asq$-linear 
graded inner products on  $\Omega(\sq)$:
\beq
\hs{\theta}{\theta^{\prime}}_{\sq}^{L}:=
\hs{\theta}{\theta^{\prime}}^{L}; \qquad\qquad\qquad
\hs{\theta}{\theta^{\prime}}_{\sq}^{R}:=
\hs{\theta}{\theta^{\prime}}^{R}
\label{prosca}
\eeq
with $\theta,\theta^{\prime}\in\,\Omega(\sq)$. The analogue result to relation \eqref{Holr} is given in the following

\begin{prop}
\label{LH}
On the exterior algebra on the sphere $\sq$ endowed with the above graded left (resp. right) inner product, there exists a unique invertible left 
$\Asq$-linear 
Hodge operator $\check{L}:\Omega^{k}(\sq)\to\Omega^{3-k}(\sq)$, (resp. a unique invertible right $\Asq$-linear Hodge operator $\check{R}$) for $k=0,1,2,3$, satisfying 
\beq\label{Les}
\int_{\check{\mu}}\theta^{*}\wedge \check{L}(\theta^{\prime})=\hs{\theta}{\theta^{\prime}}^{L}_{\sq} \, , \qquad\qquad\qquad
\int_{\check{\mu}}\theta^{*}\wedge \check{R}(\theta^{\prime})=\hs{\theta}{\theta^{\prime}}^{R}_{\sq}
\eeq
for any $\theta,\theta^{\prime}\in \Omega^{k}(\sq)$.
They can be written in terms of the sesquilinear map \eqref{p4q} as: 
\beq \label{cLo}
\begin{array}{lcl}
 \check{L}(1)=\check{\mu} \, , &
 \quad\quad & \check{L}(\check{\mu})= \{\check{\mu},\check{\mu}\}  \\
 \check{L}(\phi\,\omega_{-})=\check{m}\alpha\,\phi\,\omega_{-}\wedge\omega_{0} \, , &
 \quad \quad & \check{L}(\phi\,\omega_{-}\wedge\omega_{0})=\check{m}\,\{\omega_-\wedge\omega_{0},\omega_-\wedge\omega_0\}\phi\,\omega_{-} \, ,\\
\check{L}(\phi^{\prime}\omega_{+})=\check{m}\,q^2\alpha\,\phi^{\prime}\omega_{0}\wedge\omega_{+} \, ,  & 
\quad\quad   & \check{L}(\phi^{\prime}\omega_{0}\wedge\omega_{+})=\check{m}\,\{\omega_{+}\wedge\omega_{0},\omega_{+}\wedge\omega_{0}\}\,\phi^{\prime}\omega_{+} \, , \\
\check{L}(\omega_{0})=-\check{m}\nu\,\omega_{-}\wedge\omega_{+} \, , &
\quad\quad& \check{L}(\omega_{-}\wedge\omega_{+})=-\check{m}\,\{\omega_{-}\wedge\omega_{+},\omega_{-}\wedge\omega_{+}\}\,\omega_{0} \, 
 \end{array}
 \eeq
and
\beq
\begin{array}{lcl}
 \check{R}(1)=\check{\mu} \, , &  \quad\quad & \check{R}(\check{\mu})= \{\check{\mu},\check{\mu}\}  \\
 \check{R}(\omega_{-}\,\phi)=\check{m}q^2\alpha\,\omega_{-}\wedge\omega_{0}\,\phi \, , & \quad\quad
& \check{R}(\omega_{-}\wedge\omega_{0}\,\phi)=\check{m}\,q^2\{\omega_-\wedge\omega_{0},\omega_-\wedge\omega_0\}\,\omega_{-} \,\phi ,\\
\check{R}(\omega_{+}\,\phi^{\prime})=\check{m}\,\alpha\,\omega_{0}\wedge\omega_{+} \,\phi^{\prime} , & \quad\quad   & \check{R}(\omega_{0}\wedge\omega_{+}\,\phi^{\prime})=\check{m}\,q^{-2}\{\omega_{+}\wedge\omega_{0},\omega_{+}\wedge\omega_{0}\}\,\omega_{+} \,\phi^{\prime} , \\
\check{R}(\omega_{0})=-\check{m}\nu\,\omega_{-}\wedge\omega_{+} \, , & \quad\quad& \check{R}(\omega_{-}\wedge\omega_{+})=-\check{m}\,\{\omega_{-}\wedge\omega_{+},\omega_{-}\wedge\omega_{+}\}\,\omega_{0} \, .
 \end{array}
 \label{cRo}
 \eeq

\begin{proof}
For the rather technical proof we refer to  \cite{ale09}, where the same strategy has been adopted for the analysis of an Hodge operator on a two dimensional  exterior algebra on $\sq$.  Here we only observe that the uniqueness follows from the result in Lemma~\ref{lef}. Given two  operators $\check{L},\check{L}^{\prime}:\Omega^{k}(\sq)\to\Omega^{3-k}(\sq)$ satisfying \eqref{Les} (or equivalently $\check{R},\check{R}^{\prime}$), their difference must satisfy the relation $\int_{\check{\mu}}\theta^{\prime*}\wedge (\check{L}(\theta)-\check{L}^{\prime}(\theta))=0$ for any $\theta,\theta^{\prime}\in \Omega^{k}(\sq)$. The left-faithfulness of the integral allows one  then eventually to get  $\check{L}(\theta)=\check{L}^{\prime}(\theta)$.  
\end{proof}
\end{prop}
From \eqref{p13} and \eqref{p15-} it is $\check{\mu}=\check{m}\,\chi_{z}$, so we define
\begin{align}
 \det\check{\Gamma}:=\frac{\Gamma(\chi_z,\chi_z)}{2(1+q^2+q^{-2})} ,  
 \qquad sgn(\check{\Gamma}):=\frac{\det(\check{\Gamma})}{|\det\check{\Gamma}|}
\label{pxiv22}
\end{align}
and set 
$$
\check{m}^2\det\check{\Gamma}:=sgn(\check{\Gamma})
$$ as a definition for the scale factor $\check{m}\in\,\IR$. Clearly this choice gives $\check{L}^2(1)=\check{R}^2(1)=sgn(\check{\Gamma})$. 
Analogously to what happened for $\SU$ before, the sign in \eqref{pxiv22} turns out to be positive for the class of contractions we are considering, i.e. $sgn(\check{\Gamma})=1$.

We conclude by noticing that the Hodge operators \eqref{cLo} and  \eqref{cRo}
are diagonal, but still there is no choice for the parameters \eqref{piv22} and \eqref{pv22} of a real and hermitian contraction map such that a relation like \eqref{p8} is satisfied.

\section{Laplacian operators}
\label{se:L}
Given the Hodge operators  constructed in the previous sections, the corresponding Laplacian operators on the quantum group $\SU$,
\begin{align*}
&\Box^L_{\SU}:\ASU\to\ASU , \qquad\qquad  \Box^{L}_{\SU}(x) := - \star^L\dd\star^L\dd x , \\
&\Box^R_{\SU}:\ASU\to\ASU , \qquad\qquad  \Box^{R}_{\SU}(x) := - \star^R\dd\star^R\dd x  
\end{align*}
can be readily written in terms of the basic derivations \eqref{Lq} and \eqref{Rder}  for the first order differential calculus as  
\begin{align}
\label{elp4}
 \Box^L_{\SU} x = 
\left\{\alpha \left( L_{+}L_{-}+q^2 L_{-}L_{+}\right) + \nu\,L_{0}L_{0}+\gamma\,L_{z}L_{z}+2\epsilon L_{0}L_{z}\right\}\lt x, 
\end{align}
and
\begin{align}
 \Box^R_{\SU} x = 
\left\{\alpha \left( q^2 R_{+}R_{-}+R_{-}R_{+} \right) + \nu\,R_{0}R_{0}+\gamma\,R_{z}R_{z}+2\epsilon R_{0}R_{z}\right\}\lt x, 
\label{erp4}
\end{align}
with parameters given in \eqref{piv22}.

From the decomposition \eqref{dcmp} and the action \eqref{rellb} it is immediate to see that such Laplacians restrict to  operators $:\cl_{n}\to\cl_{n}$. In order to diagonalise them, we recall the decomposition \eqref{decoln}.   
The action of each term of the Laplacians on the  basis elements $\{\phi_{n,J,l} \}$ in \eqref{bsphi} can be explicitly computed by \eqref{lact}, giving:
\begin{align}
&L_{-}L_{+}\lt\,\phi_{n,J,l}=q^{-1-n} \, [J-\half n][J+1- \half n ]\,\phi_{n,J,l} \, , \nn \\
&L_{+}L_{-}\lt\,\phi_{n,J,l}=q^{1-n} \, [J+\half n][J+1-\half n]\,\phi_{n,J,l} \, , \nn \\
&L_{z}\lt\,\phi_{n,J,l}=-q^{-\frac{1}{2} n } \, [\half n]\,\phi_{n,J,l} \, , \nn \\
&L_{0}\lt\,\phi_{n,J,l}=([J+\half]^{2}-[\half]^2)\,\phi_{n,J,l}=[J][J+1]\,\phi_{n,J,l} \, .
\label{eopep}
\end{align} 
Here for the labels one has $n\in\IN$ with $J=\tfrac{|n|}{2}+\IZ$ and $l=0,\ldots,2J$. 

The Laplacians on the quantum sphere are, with  $f\in\,\Asq$: 
\begin{align}
\Box_{\sq}^L f:=-\check{L}\dd\check{L}\dd f= 
\left\{\alpha\,L_{+}L_{-}+q^2\alpha\,L_{-}L_{+}+\nu\,L_{0}L_{0}\right\}\lt f, 
\label{sfl} 
\end{align}
and
\begin{align}
\Box_{\sq}^R f:=-\check{R}\dd\check{R}\dd f= 
\left\{q^2 \alpha\,R_{+}R_{-}+\alpha\,R_{-}R_{+}+\nu\,R_{0}R_{0}\right\}\lt f. 
\label{sfr} 
\end{align}
They both are the restriction to $\sq$ of the Laplacian on $\SU$, the left and right one respectively.
Their actions can be written in terms of the action of the Casimir element $C_{q}$ of $\su$, immediately giving their spectra. They coincide on $\sq$:
\begin{align*}
\Box_{\sq}^{L,R}  &=  2q\alpha(C_{q}+\tfrac{1}{4}-[\tfrac{1}{2}]^2)
+\nu(C_{q}+\tfrac{1}{4}-[\tfrac{1}{2}]^2)^2 ,  \nn \\
&=  2q\alpha \,L_{0} +\nu  \, L^2_{0} \qquad \mathrm{on} \quad \Asq . 
\end{align*}
Using \eqref{eopep},  spectra are readily found:
\beq
 \Box_{\sq}^{L,R} (\phi_{0,J,l}) = \left( 2q\alpha[J][J+1]+\nu[J]^2[J+1]^2 \right) \phi_{0,J,l} ,
\label{specL}
\eeq
with $J\in\,\IN, l=0,\ldots,2J$. 
We end this section by comparing these spectra to the spectrum of $D^2$, the square of the Dirac operator on $\sq$ studied in \cite{gs10}. Some straightforward computation leads to:
\begin{align}\label{cosp}
\text{spec}( \,\Box_{\sq}^{L,R})=\text{spec}(D^2-[\tfrac{1}{2}]^2)\quad\quad\Leftrightarrow\quad\quad2q\alpha=1, \quad \nu=q^{-2}(q-q^{-1})^4 .
\end{align}


\section{A digression: connections  on the Hopf fibration over $\sq$ } \label{s:ccd}


A monopole connection for the quantum fibration $\Asq\hookrightarrow\ASU$ on the standard Podle\'s sphere -- with a left-covariant 3d calculus on $\SU$ and the (corresponding restriction to a) 2d left-covariant calculus on $\sq$ -- was explicitly described in \cite{bm93}. 
A slightly different, but to large extent equivalent \cite{durcomm} formulation of this and of a fibration constructed on the same topological data $\Asq\hookrightarrow\ASU$, but with $\SU$ equipped with a bicovariant 4D calculus inducing on 
$\sq$ a left-covariant 3d calculus, are presented in \cite{durII}.  
The general problem of finding the conditions between the differential calculi on a base space algebra and on a `structure' group, in a way giving a principal bundle structure with compatible calculi and a consistent definition of connections on it has been deeply studied \cite{bm98,ps97,dur98,haj97}. The slightly different perspective of this digression is to follow  the path reviewed in appendix~\S\ref{ap:qpb},  namely to recall from \cite{gs10} the formulation  of a Hopf bundle on the standard Podle\'s sphere starting from the 4D bicovariant calculus \`a la Woronowicz on the total space $\SU$, in order to fully describe  the set of its connections.  The first step in this analysis consists in describing how the differential calculus on $\SU$ naturally induces a 1 dimensional bicovariant calculus on the structure group $\U(1)$, and in which sense these two calculi are compatible.

\subsection{A 1D bicovariant calculus on $\U(1)$}
\label{se:cu1}

The Hopf projection \eqref{qprp} allows one to define an ideal $\cq_{\U(1)}\subset\ker\varepsilon_{\U(1)}$ as the projection
$\cq_{\U(1)}=\pi(\cq_{\SU})$. Then $\cq_{\U(1)}$ is generated by the three elements 
\begin{align*}
\xi_{1}&=(z^{2}-1)+q^{2}(z^{-2}-1), \\ 
\xi_{2}&=(q^{2}z+z^{-1}-(q^{3}+q^{-1}))(q^{2}z+z^{-1}-(1+q^{2})), \\ 
\xi_{3}&=(q^{2}z+z^{-1}-(q^{-1}+q^{3}))(z^{-1}-z), 
\end{align*}
and, since $\mathrm{Ad}(\cq_{\U(1)})\subset\cq_{\U(1)}\otimes\ca(\U(1))$, it corresponds to a bicovariant differential calculus on $\U(1)$. The identity
$$
-q(1+q^{4})^{-1}(1+q^{2}+q^{3}+q^{5})^{-1}\{(q^{6}-1)\xi_{3}+(1+q^{4})\xi_{2}-q^{2}(1+q^{2})\xi_{1}\}=(z-1)+q(z^{-1}-1)
$$
shows that $\xi=(z-1)+q(z^{-1}-1)$ is in $\cq_{\U(1)}$. By induction one also sees that  
\begin{align}
&j>0:\qquad z^{j}(z-1)=\xi(\sum_{n=0}^{j-1}\,q^{n}z^{j-n})+q^{j}(z-1), \nn \\
&j<0:\qquad z^{-\mj}(z-1)=-\xi(\sum_{n=1}^{\mj-1}\,q^{-n}z^{n-\mj})+q^{-\mj}(z-1).
\label{idzj}
\end{align}
From these relations it is immediate to prove (as in \cite{gs10}) that there is a complex vector space isomorphism 
$\ker\varepsilon_{\U(1)}/\cq_{\U(1)} \simeq\IC$.  The differential calculus induced by $\cq_{\U(1)}$ is 1-dimensional, and the projection $\pi_{\cq_{\U(1)}}:\ker\varepsilon_{\U(1)}\to\ker\varepsilon_{\U(1)}/\cq_{\U(1)}$ can be written as
\beq
\label{lam4}
\pi_{\cq_{\U(1)}}:\quad z^{j}(z-1)\to q^{j}[z-1],
\eeq
on the vector space basis $\varphi(j)=z^{j}(z-1)$ in $\ker\varepsilon_{\U(1)}$, with notation $[z-1]\,\in\,\ker\varepsilon_{\U(1)}/\cq_{\U(1)}$. The projection \eqref{lam4} will be used later on to define connection 1-forms on the fibration.

As a basis element for the quantum tangent space $\mathcal{X}_{\cq_{\U(1)}}$ we take 
\beq
X=L_z=\frac{K^{-2}-1}{q-q^{-1}}.
\label{X1}
\eeq
The $*$-Hopf algebras $\ca(\U(1))$ and $\cu(1)\simeq\{K,K^{-1}\}$ are dually paired via the pairing, induced by the one in \eqref{ndp} between $\ASU$ and $\su$, with 
\beq
\hs{K^{\pm1}}{z}=q^{\mp\frac{1}{2}},\qquad\hs{K^{\pm1}}{z^{-1}}=q^{\pm\frac{1}{2}},
\label{bp1}
\eeq
on the generators. Thus, the exterior derivative $\dd:\ca(\U(1))\to\Omega^{1}(\U(1))$ can be written, 
for any $u\in\,\ca(\U(1))$,  
as $ \dd u=(X\lt u)\ \theta $ on the left invariant basis 1-form $\theta\sim[z-1]$. On the generators of the coordinate algebra one has 
\beq
\dd z=\frac{q-1}{q-q^{-1}}\,z \ \theta, 
\qquad \dd z^{-1}=\frac{q^{-1}-1}{q-q^{-1}}\,z^{-1} \ \theta,
\label{dz1}
\eeq
so to have $\theta=(q-1)(q-q^{-1})^{-1}z^{-1}\dd z$. From the coproduct $\Delta X=1\otimes X+X\otimes K^{-2}$ the 
$\ca(\U(1))$-bimodule structure in $\Omega^{1}(\U(1))$ is 
$$
\theta \ z^{\pm}=q^{\pm}z^{\pm}\ \theta.
$$ 

\subsection{Connections on the principal bundle}\label{conHf}
The compatibility -- as described in appendix~\S\ref{ap:qpb} and expressed by the exactness of the sequence \eqref{des} -- of the differential calculus $\U(1)$ presented above with the 4D differential calculus on $\SU$ presented in \S\ref{se:4dc}, has been proved in \cite{gs10}. As a consequence, collecting the various terms, the data  
$$
\left(\ASU, \Asq, \ca(\U(1)); \cn_{\SU}=r^{-1}(\SU\otimes\cq_{\SU}), \cq_{\U(1)}\right)
$$ 
is a quantum principal bundle with the described calculi. 

In order to obtain connections on this bundle, that is maps \eqref{si} splitting the sequence \eqref{des}, we need to compute the action of the map 
$$\sim_{\cn_{\SU}}:\Omega^{1}(\SU)\to\ASU\otimes(\ker\varepsilon_{\U(1)}/\cq_{\U(1)})$$ defined via the diagram \eqref{qdia}.
Since it is left $\ASU$-linear, we take as  representative universal 1-forms corresponding to  the left invariant 1-forms  \eqref{om4} in $\Omega^1(\SU)$:
\begin{align*}
&\pi^{-1}_{\cn_{\SU}}(\omega_{+})=(a\delta c-qc\delta a) \\ 
&\pi^{-1}_{\cn_{\SU}}(\omega_{-})= (c^* \delta a^*- qa^*\delta c^*) \\ 
&\pi^{-1}_{\cn_{\SU}}(\omega_{0})=\{a^{*}\delta a+c^{*}\delta c+ q(a\delta a^{*}+q^{2}c\delta c^{*})\}/ (q+1)\lambda_{1}  \\
&\pi^{-1}_{\cn_{\SU}}(\omega_{z})=a^{*}\delta a+c^*\delta c - (a\delta a^*+q^2 c\delta c^*) .
\end{align*}
On them the  action of the canonical map \eqref{chimap} is found to be:
\begin{align*}
&\chi(a\delta c-qc\delta a)=(ac-qca)\otimes(z-1)=0 \\ 
&\chi(c^* \delta a^*- qa^*\delta c^*) = (c^*a^* - qa^*c^*)\otimes(z^{*}-1)=0 \\
&\chi\left((1+q)^{-1}\lambda_{1}^{-1}\{a^{*}\delta a+c^{*}\delta c+ q(a\delta a^{*}+q^{2}c\delta c^{*})\}\right)=1\otimes\{(z-1)+q(z^{-1}-1)\}=1\otimes\xi \\
&\chi(a^{*}\delta a+c^*\delta c - (a\delta a^*+q^2 c\delta c^*))=1\otimes(z-z^{-1})
\end{align*}
with $\xi\in\cq_{\SU}$ introduced in \S\ref{se:cu1}. From the  isomorphism \eqref{lam4} one finally has: 
\begin{align}
&\sim_{\cn_{\SU}}(\omega_{\pm})= \sim_{\cn_{\SU}}(\omega_{0})=0 \nn \\
&\sim_{\cn_{\SU}}(\omega_{z})=1\otimes(1+q^{-1})[z-1].
\label{omzv}
\end{align} 
{}From these one recovers $\Omega^{1}_{\mathrm{hor}}(\SU)=\ker\sim_{\cn_{\SU}}$ with, using \eqref{bis2}, 
\beq
\label{le:hf}
\ker\sim_{\cn_{\SU}}\simeq\ASU\{\omega_{\pm},\omega_{0}\}\simeq\{\omega_{\pm},\omega_{0}\}\ASU.
\eeq
 
\begin{rema} \label{horv}
\textup{
From \eqref{omzv}, for the generator $X=L_z$  in \eqref{X1} one gets that 
$$
\widetilde{X}(\omega_z) = \langle{X},{\sim_{\cn_{\SU}}(\omega_{z})}\rangle
=1 ,
$$ 
which identifies $L_{z}\in\,\mathcal{X}_{\cq}$ as a vertical vector for the fibration. In turn it is used to extend the notion of horizontality to higher order forms in $\Omega(\SU)$. One defines \cite{KS97} a contraction operator $i_{L_{z}}:\Omega^{k}(\SU)\to\Omega^{k-1}(\SU)$, giving  $i_{L_{z}}(\omega_{\pm})=i_{L_{z}}(\omega_{0})=0$, and $i_{L_{z}}(\omega_{z})=1$ on 1-forms, so  that $\ker i_{L_{z}}\simeq\Omega^{1}_{\mathrm{hor}}(\SU)$. Then one defines
\beq
\Omega^{k}_{\mathrm{hor}}(\SU):=\left.\ker\right|_{\Omega^{k}(\SU)}i_{L_{z}},
\label{khor}
\eeq
 that is the kernel of the contraction map when restricted to the bimodule of $k$-forms. 
}
\end{rema}

Given the explicit expression \eqref{omzv} for the canonical map compatible with the differential calculi we are using, and the $\ca(\U(1))$-coaction  
\beq
\delta_{R}^{(1)}\omega_{z}=\omega_{z}\otimes 1, \qquad
\delta_{R}^{(1)}\omega_{0}=\omega_{0}\otimes 1, \qquad
\delta_{R}^{(1)}\omega_{\pm}=\omega_{\pm}\otimes z^{\pm2},
\label{rifo}
\eeq 
using the vector space basis $\varphi(j)$ in $\ker\varepsilon_{\U(1)}$ of \S\ref{se:cu1}, 
a connection \eqref{si} is given by 
\beq
\sigmat(\phi\otimes[\varphi(j)])=q^{-2j}(1+q^{-1})^{-1}\phi(\omega_{z}+\mathrm{a})
\label{si3}
\eeq
for any $\phi\in\,\ASU$ and any element $\mathrm{a}\in\,\Omega^1(\sq)$. On vertical forms, the projection $\Pi$ associated to this connection turns out  to be
\begin{align}
&\Pi(\omega_{\pm})=0=\Pi(\omega_{0}),\nn \\
&\Pi(\omega_{z})=\sigmat(\sim_{\cn_{\SU}}(\omega_{z}))=\sigmat(1\otimes[\varphi(0)])=\omega_{z}+\mathrm{a},
\label{Pi3}
\end{align}
while the corresponding connection 1-form $\omega:\ca(\U(1))\to\Omega^{1}(\SU)$ is given by 
\beq
\omega(z^{n})=\sigmat(1\otimes[z^{n}-1]) 
=q^{n/2}[\tfrac{n}{2}] (\omega_{z}+\mathrm{a}).
\label{ome3}
\eeq
Connections corresponding to $\mathrm{a}=s\omega_{0}$ with $s\in\,\IR$ were already considered in \cite{durII}.

The vertical projector \eqref{Pi3}  allows one to define a covariant derivative 
$$
\mathfrak{D}:\ASU\to\Omega^{1}_{\mathrm{hor}}(\SU), 
$$
given (as usual) as the horizontal projection of the exterior derivative:
\beq
\label{mfD}
\mathfrak{D}\phi:=(1-\Pi)\dd\phi.
\eeq
Covariance here clearly refers to  the right  coaction of the structure group $\U(1)$ of the bundle,  since it is that $\delta_{R}\phi=\phi\otimes z^{-n} \,\Leftrightarrow\,\delta_{R}^{(1)}(\mathfrak{D}\phi)=(\mathfrak{D}\phi)\otimes z^{-n}$. From \eqref{ome} the action of this operator can be written as 
\beq
\label{mfDo}
\mathfrak{D}\phi=\dd\phi-\phi\wedge\omega(z^{-n})
\eeq
 for any $\phi\in\,\cl_{n}$. From the bimodule structure \eqref{bis2} it is easy to check that all the above connections are \emph{strong} connections  in the sense of \cite{hajsc}.

\medskip

The analysis in this  section allows us to prove the results in Proposition~\ref{lea}. The exterior algebra $\Omega(\sq)$ is defined to be the set of horizontal and $\U(1)$-coinvariant elements in $\Omega(\SU)$, with respect to the extension $\delta_{R}^{(k)}$ (introduced in App.\ref{ap:qpb}) of the canonical coaction \eqref{cancoa} to higher order forms in $\Omega(\SU)$. It is then easy to check, from \eqref{khor} and \eqref{rifo}, that the isomorphisms given in  expressions  \eqref{ntd} for $\Omega(\sq)$ do hold.

\section{Gauged Laplacians on line bundles} \label{GLoLB}





Each $\Asq$-bimodule $\cl_{n}$ defined in \eqref{libu} is a bimodule of co-equivariant elements in $\ASU$ for the right $\U(1)$-coaction \eqref{cancoa}, and as such  can be thought of as a module of `sections of a line bundle' over the quantum sphere $\sq$. Without requiring any compatibility with additional structures, any $\cl_{n}$ can be realized both as a projective right or left $\Asq$-module (of rank 1 and winding number $-n$). One of such structures is that of a connection
on the quantum principal bundle $\Asq\hookrightarrow\ASU$. By transporting the covariant derivative \eqref{mfD} on the principal bundle to a derivative on sections forces to break the symmetry between  the left or the right  $\Asq$-module realization of $\cl_{n}$. 

With the choice in 
\S\ref{s:qsb} for the principal bundle, we need an isomorphism  $\cl_{n}\simeq\cf_n$ with $\cf_{n}$ a projective left $\Asq$-module \cite{HM98}. 
This isomorphism is constructed in terms of a projection operator $\qpp$.
Given this identification, in \S\ref{se:lbss} we shall describe the complete equivalence between covariant derivatives on $\mathcal{F}_{n}$ (associated  to the $3\dd$ left covariant differential calculus over 
$\sq$)  and connections (as described in \S\ref{s:ccd}) on the  principal bundle $\Asq\hookrightarrow\ASU$, corresponding to  compatible   $4D_{+}$ bicovariant calculus over $\SU$ and  $3\dd$ left covariant calculus over $\sq$. We shall then move  to  a family of gauged Laplacian operators on $\mathcal{F}_{n}$, obtained by coupling the Laplacian operator over the quantum sphere with a set of suitable gauge potentials. We
finally show that among them there is one whose action extends to $\cl_n$ the action of the Laplacian \eqref{sfr} on $\cl_0\simeq\sq$. As we noticed in \S\ref{se:L}, the action of the (right) Laplacian \eqref{sfr} on $\sq$ is given by the restriction of the action \eqref{erp4} of the (right) Laplacian $\Box^R_{\SU}$.  Here we obtain that the action of such gauged Laplacian can be written in terms of the ungauged (right) Laplacian on $\SU$, in parallel to what happens on a classical principal bundle  (see e.g. \cite[Prop.~5.6]{bgv}) and on the Hopf fibration of the sphere $\sq$ with calculi coming from the left covariant one on $\SU$ as shown in \cite{lareza,ale09}. 

%

\subsection{Line bundles as projective left $\Asq$-modules}\label{se:lbss}
Every (equivariant, the only ones we use in this paper) finitely generated projective (left or right) $\Asq$-module is a direct sum of  $\cl_{n}$'s (cf. \cite{SWPod}). As said, these are line bundles of degree $-n$ on the quantum sphere and to describe them all one needs is a collection of idempotents $\qpp$, which we are going to introduce.
With $n\in\,\IZ$, we consider the projective left
$\Asq$-module $\cf_{n}=(\Asq)^{\mn+1}\qpp$, with projections  \cite{BM98,HM98} (cf. also \cite{lareza})
\beq
\qpp=\ket{\Psin}\bra{\Psin} 
\label{dP},
\eeq
written in terms of elements 
$\ket{\Psin}\in \ASU^{\mn+1}$ and their duals $\bra{\Psin}$ as follows. One has: 
\begin{align}
&n\leq0:\qquad\ket{\Psin}_{\mu}=\sqrt{\alpha_{n,\mu}}  ~ c^{\mn-\mu}a^{\mu}\,\in \cl_{n}, \nn \\
&\mathrm{where}\qquad\alpha_{n,\mn}=1;\qquad\alpha_{n,\mu}=\prod\nolimits_{j=0}^{\mn-\mu-1}\left(\frac{1-q^{2\left(\mn-j\right)}}
{1-q^{2\left(j+1\right)}}\right),  \quad \mu = 0, \ldots,  \mn-1  
\label{ketpi}
\end{align}
\begin{align}
&n\geq 0:\qquad\ket{\Psin}_{\mu}=\sqrt{\beta_{n,\mu}} ~ c^{* \mu}a^{* n-\mu}\,\in \cl_{n}, \nn\\
&\mathrm{where}  \qquad \beta_{n,0}=1;\qquad
\beta_{n,\mu}=q^{2\mu}\prod\nolimits_{j=0}^{\mu-1}\left(\frac{1-q^{-2\left(n-j\right)}}{1-q^{-2\left(j+1\right)}}\right),
\quad \mu = 1, \ldots,  n 
\label{ketni} .
\end{align}
The coefficients are chosen so that $\hs{\Psin}{\Psin}=1$, as a consequence $(\qpp)^{2}=\qpp$. Also by construction it holds that $(\qpp)^{\dagger}=\qpp$.  

The isomorphism $\cl_n \simeq \cf_{n}=(\Asq)^{\mn+1}\qpp$ is realized as follows:  
\beq
\cl_{n} ~\xrightarrow{~\simeq~}~ \cf_{n}, \quad \phi  \mapsto \bra{\sigma_{\phi}} = \phi \bra{\Psi^{(n)} },
\label{iso1n}
\eeq 
with inverse
$$
\cf_{n} ~\xrightarrow{~\simeq~}~ \cl_{n}, \quad \bra{\sigma_{\phi}}  \mapsto \phi =\hs{\sigma_{\phi}}{ \Psi^{(n)} }.
$$
Given the exterior algebra $(\Omega(\sq),\dd)$ on the quantum sphere we are considering, 
a covariant derivative on the left $\Asq$-modules  $\cf_n$ is a $\IC$-linear map
\beq
\nabla:\Omega^k(\sq)\otimes_{\Asq}\cf_{n}\,\to\,\Omega^{k+1}(\sq)\otimes_{\Asq}\cf_{n}
\label{cdev}
\eeq
that satisfies the \emph{left} Leibniz rule 
$$
\nabla(\xi\wedge\bra{\sigma})=(\dd\xi)\wedge\bra{\sigma}+(-1)^m \xi\wedge \nabla\bra{\sigma}
$$
for any $\xi\,\in\,\Omega^{m}(\sq)$ and $\bra{\sigma}\,\in\,\Omega^{k}(\sq)\otimes_{\Asq}\cf_{n}$.
The curvature associated to a covariant derivative is 
$\nabla^{2}:\cf_{n}\to\Omega^{2}(\sq)\otimes_{\Asq}\cf_{n}$, that is $\nabla^{2}(\xi\,\bra{\sigma})=\xi\,\nabla^2(\bra{\sigma})=\xi\,F_{\nabla}(\bra{\sigma})$ with the last equality defining the curvature 2-form $F_{\nabla}\in\,\textup{Hom}_{\Asq}(\cf_{n},\Omega^{2}(\sq)\otimes_{\Asq}\cf_{n})$. 

Any covariant derivative -- an element in $C(\cf_{n})$ -- and its  curvature
can be written as
\begin{align}
\nabla\bra{\sigma}&=(\dd\bra{\sigma})\qpp\,+(-1)^k\,\bra{\sigma}\An, \nn 
\\
\nabla^{2}\bra{\sigma}&=\bra{\sigma}\{-\dd\qpp\wedge\dd\qpp\,+\,\dd\An\,-\,\An\wedge\An\}\qpp .
\label{cocug}
\end{align}
with $\bra{\sigma}\in\,\Omega^{k}(\sq)\otimes_{\Asq}\cf_{n}$. The negative signs in the second expression above come from the left Leibniz rule, since form valued sections  are elements of projective left $\Asq$-modules.
For the `gauge potential' $\An$ one has 
\beq
\An=\qpp \An=\An\qpp=\ket{\Psin}\an\bra{\Psin}\,\in\,\mathrm{Hom}_{\Asq}(\cf_n, \Omega^{1}(\sq)\otimes_{\Asq}\cf_{n}),
\label{apA}
\eeq
with  $\an\in\,\Omega^{1}(\sq)$. The monopole (Grassmann) connection corresponds to $\an=0$.

In analogy with the identification \eqref{iso1n}, the covariant derivative $\nabla$ naturally induces an  operator 
$D:\cl_{n}\to\cl_{n}\otimes_{\Asq}\Omega^{1}(\sq)$  that can be written as
\beq
D\phi:=(\nabla\bra{\sigma_\phi})\ket{\Psin}=\dd\phi-\phi\,\{\hs{\Psin}{\dd\Psin}-\an\}.
\label{domeg}
\eeq 
We refer to the 1-form  
\beq
\label{coAp}
\Omega^{1}(\SU)\,\ni\,\varpi^{(n)}=\left(\hs{\Psin}{\dd\Psin}-\an\right)
\eeq
 as the  connection 1-form of the gauge potential. It   allows   to express the curvature as
\beq
F_{\nabla}=-\ket{\Psin}\left( \dd\varpi^{(n)}+\varpi^{(n)}\wedge\varpi^{(n)} \right)\bra{\Psin} 
\label{Fp}
\eeq
where  $(\dd\varpi^{(n)}+\varpi^{(n)}\wedge\varpi^{(n)})\in\,\Omega^{2}(\sq)$.

The covariant derivatives defined above on the left modules $\cf_n$ fit in the general theory of connections on the quantum Hopf bundle as described in the \S\ref{conHf}: any covariant  vertical projector, as in \eqref{Pi3}, induces a  gauge potential $\An$ as in \eqref{apA}.  The notion \eqref{coAp} of connection 1-form of a given gauge potential in $C(\cf_{n})$ matches the notion \eqref{ome3} of connection 1-form $\omega:\ca(\U(1))\to\Omega^{1}(\SU)$ on the Hopf bundle. From the $\Asq$-bimodule isomorphism $\oplus_{n\in\,\IZ}\cl_{n}\otimes_{\Asq}\Omega^{1}(\sq)\,\simeq\,\Omega^{1}_{\mathrm{hor}}(\SU)$ (see Remark~\ref{horv}), this matching amounts to equate the actions of the covariant derivative operators \eqref{domeg} and \eqref{mfD},
\beq
\forall\,\phi\in\,\cl_{n}: \quad D\phi=\mathfrak{D}\phi\qquad\Leftrightarrow\qquad\varpi^{(n)}=\omega(z^{-n}).
\label{omom}
\eeq
From formula \eqref{ome3}, this correspondence can be written as
\beq
\an=\lambda_{n}\omega_{0}-\xi_{-n}\mathrm{a},
\label{figi}
\eeq
where the coefficients refer to the eigenvalue equations:
\begin{align}
&L_{z}\lt\ket{\Psin}:=\xi_{n}\ket{\Psin} \qquad\Rightarrow\qquad\xi_{n}=-q^{-\tfrac{n}{2}}\,[\frac{n}{2}]\nn \\
&L_{0}\lt\ket{\Psin}:=\lambda_{n}\ket{\Psin}\qquad\Rightarrow\qquad\lambda_{n}=[\frac{\mn}{2}] \label{lamni}
[\frac{\mn}{2}+1].
\end{align}
Finally, the equivalence \eqref{omom} allows one to  introduce a covariant derivative 
$$ D:\Omega^{k}_{\mathrm{hor}}(\SU)\to\Omega^{k+1}_{\mathrm{hor}}(\SU),
$$ thus extending to horizontal  forms on the total space of the quantum Hopf bundle the covariant derivative operator on $\ASU$ as given in \eqref{mfD}. This  follows the formulation described in \cite{hajsc}, since any connection on the principal bundle is strong. Upon defining 
$$
\cl_{n}^{(k)}:=\{\phi\in\,\Omega^{k}_{\mathrm{hor}}(\SU)\,:\,\delta_{R}^{(k)}\phi=\phi\otimes z^{-n}\},
$$
where $\delta_{R}^{(k)}$ is the natural right  $\U(1)$-coaction on $\Omega^{k}(\SU)$, one obtains:
\beq
D\phi=\dd\phi-(-1)^{k}\phi\wedge\omega(z^{-n}).
\label{esD}
\eeq 
A further  extension  to the whole exterior algebra $\Omega(\SU)$ is proposed in \cite{durII}: a generalisation of the analysis in \cite[\S9]{ale09} shows how this extension is far from being unique. 

We restrict  our analysis again to covariant derivatives $\nabla_{s}\ket{\sigma}$ in \eqref{cocug}
whose gauge potential and corresponding connection 1-form 
are of the form:
\begin{align}
\Ans=s\ket{\Psin}\omega_{0}\bra{\Psin},  \qquad\qquad \varpi_{s}^{(n)}=\xi_{n}\omega_{z}\,+\,(\lambda_{n}-s)\omega_{0} ,
\label{osn}
\end{align}
for $s\in\IR$ and coefficients as in \eqref{lamni}, since they reduce in the classical limit to the  monopole connection on line bundles associated to the classical Hopf bundle $\pi:S^3\to S^{2}$.
Relations \eqref{2fw} and \eqref{d2f} allow to compute the curvature 2-form \eqref{Fp} as
\begin{align*}
&\dd\varpi_{s}^{(n)}=\left( (q+q^{-1})\xi_{n}-(s-\lambda_{n})(q-q^{-1}) \right) \omega_+\wedge\omega_- \nn \\
&\varpi_{s}^{(n)}\wedge\varpi_{s}^{(n)}=(q-q^{-1})\xi_{n} \left( (q+q^{-1})\xi_{n}+(q-q^{-1})(s-\lambda_{n}) \right)
\omega_{+}\wedge\omega_{-}.
\end{align*}

\subsection{Gauged Laplacians}

In order to introduce an Hodge operator 
\beq
\star^{\cR}:\Omega^{k}(\sq)\otimes_{\Asq}\cf_{n}\to\Omega^{3-k}(\sq)\otimes_{\Asq}\cf_{n},
\label{scR}
\eeq
 we use the right $\Asq$-linear Hodge operator \eqref{cRo} on $\Omega(\sq)$:
\beq
\star^{\cR}(\xi\,\bra{\sigma}):=(\check{R}\xi)\bra{\sigma}
\label{scRe}
\eeq
 so that  a gauged Laplacian operator is  defined as 
 $$
 \Box^{\cR}_{\nabla}:\cf_{n}\to\cf_{n},\qquad\qquad\qquad\Box^{\cR}_{\nabla}\bra{\sigma}:=-\star^{\cR}\nabla(\star^{\cR}\nabla\bra{\sigma}) .
$$ 
Equivalently we have an  operator  on  $\cl_{n}\simeq\cf_{n}$ via the left $\Asq$-modules isomorphism \eqref{iso1n}. With $\phi=\hs{\sigma}{\Psin}$,  it holds that 
\beq
\Box^{\cR}_{\nabla}:\cl_{n}\to\cl_{n}, \qquad\qquad\qquad \Box^{\cR}_{\nabla}\phi=(\Box^{\cR}_{\nabla}\bra{\sigma})\ket{\Psin}.
\label{glp}
\eeq
With the family of connections \eqref{osn} and using the identities
\begin{align}
&\left(R_{\pm}\lt\,\bra{\sigma}\right)\ket{\Psin}=q^{-n}\,R_{\pm}\lt \phi ,\nn \\
&\left(R_{0}\lt\,\bra{\sigma}\right)\ket{\Psin}=q^{-n}\left(R_{0}-[\frac{\mn}{2}][1-\frac{\mn}{2}]\right)\lt\phi
\label{idL}
\end{align}
one readily computes:
\beq
\label{rfB}
\Box_{\nabla_{s}}^{\cR}\phi=
q^{-2n}\left\{ \alpha\, ( q^2 R_{+}R_{-}+ R_{-}R_{+} ) +\nu\,(R_{0}+sq^{-n}-[\frac{\mn}{2}][1-\frac{\mn}{2}])^2\right\}\lt \phi.
\eeq
Finally, fixing the parameter to be 
\beq
\label{sCr}
s(n)= q^{n} [\frac{\mn}{2}][1-\frac{\mn}{2}], 
\eeq
the action of the gauged Laplacians   extends, apart from a multiplicative factor depending on the label $n$,  to elements in the line bundles  $\cl_{n}$ the action of the Laplacian operator  \eqref{sfr}  on the quantum sphere, that is, 
\beq
q^{2n}\left(\Box^{\cR}_{\nabla_{s}}\phi\right)
=\left\{ \alpha ( q^2 R_{+}R_{-}\,+\,R_{-}R_{+} )\,+\,\nu R^2_{0} \right\}\lt\phi.
\label{finr}
\eeq
From \eqref{Rder}, the above action can be written on $\phi\in\,\cl_{n}$ as the left action \eqref{d4} of a polynomial in $\su$.  We get
\begin{align}\label{zed}
\Box^{\cR}_{\nabla_{s}} &= (- 2q \alpha\, R_0 K^{2}+\nu\,R_{0}^2) K^{-4}  - q\alpha
\frac{(q+q^{-1})(K-K^{-1})^2}{(q-q^{-1})^2}\,K^{-2}  \\
&=2q \alpha\, (C_{q} - [\tfrac{1}{2}]^{2}+\tfrac{1}{4} ) K^{-2}+\nu\,(C_{q} - [\tfrac{1}{2}]^{2}+\tfrac{1}{4} )^2) K^{-4}  - q\alpha
\frac{(q+q^{-1})(K-K^{-1})^2}{(q-q^{-1})^2}\,K^{-2} \nn ,
\end{align}
having used \eqref{casbis}. This relation is the counterpart of what happens on a classical principal bundle  (see e.g. \cite[Prop.~5.6]{bgv}) and on the Hopf fibration of the sphere $\sq$ with calculi coming from the left covariant one on $\SU$ as shown in \cite{lareza,ale09}.

\appendix

\section{Exterior differential calculi  on Hopf algebras}\label{ass:a1}

In this appendix we briefly recall general definitions and results from the theory of differential calculi on quantum spaces and quantum groups. We confine to the notions that we need in this paper, in order to construct  the exterior algebras over the quantum group $\SU$ and its subalgebra $\sq$. For a more complete analysis we refer to \cite{wor89,KS97}.

Let $\ca$ be a unital $*$-algebra over $\IC$ and $\Omega^{1}(\ca)$ an $\ca$-bimodule. Given the  linear map $\dd:\ca\to\Omega^1(\ca)$, the pair $(\Omega^1(\ca), \dd)$ is a (first order) differential calculus over $\ca$  if $\dd$ satisfies the Leibniz rule, $\dd (x y) = (\dd  x) y  + x \dd y $ for $x,y\in \ca$, and 
if $\Omega^1(\ca)$ is generated by $\dd(\ca)$ as a $\ca$-bimodule. Furthermore, it is a $*$-calculus if there is an anti-linear involution $*:\Omega^1(\ca)\to\Omega^1(\ca)$ such that $(a_{1}(\dd a)a_{2})^*=a_{2}^{*}(\dd(a^*))a_{1}^{*}$ for any $a,a_{1},a_{2}\in \ca$. 

The universal calculus $(\Omega^1(\ca)_{un}, \delta)$ 
has universal 1-forms given by the submodule $\oca{1}_{un} = \ker(m:\ca \otimes \ca \to \ca) \subset  \ca \otimes \ca$ with $m(a\otimes b)=ab$ the multiplication map, while  the universal differential 
$\delta: \ca \to \oca{1}_{un}$ is $\delta a = 1 \otimes a - a \otimes 1$. It is universal since given
any sub-bimodule $\cn$  of $\oca{1}_{un}$ with
projection $\pi_{\cn}: \oca{1}_{un} \to \oca{1} = \oca{1}_{un} /\cn$, then $(\oca{1}, \dd)$, with $\dd:=\pi_{\cn}\circ\delta$, is a first order
differential calculus over $\ca$ and any such a calculus can be
obtained in this way. The projection $\pi_{\cn}:\Omega^{1}(\ca)_{un}\to\Omega^{1}(\ca)$ is $\pi_{\cn}(\sum_{i}a_{i}\otimes b_{i})=\sum_{i}a_{i}\dd b_{i}$ with associated subbimodule $\cn=\ker\pi$.

Next, suppose $\ca$ is a left $\ch$-comodule algebra for the quantum group $\ch=(\ch, \Delta, \varepsilon, S)$,  with left  coaction $\delta_L : \ca \to \ch \otimes \ca$, an algebra map. 
The calculus is said to be left covariant provided a left coaction $\delta^{(1)}_{L}:
\oca{1} \to \ch\otimes\oca{1}$ exists, such that $\delta^{(1)}_{L}(\dd a)=(1 \otimes \dd) \delta_{L}(a)$ and $\delta^{(1)}_{L}(a_{1}\alpha\,a_{2})=\delta_{L}(a_{1})\,\delta^{(1)}_{L}(\alpha)\,\delta_{L}(a_{2})$ for any $\alpha\in\,\Omega^{1}(\ca)$ and $a_{1},a_{2}\in\,\ca$.

Left covariance of a calculus can be stated in terms of the subbimodule $\cn\subset\Omega^1(\ca)$. The left  coaction $\delta_{L}$ is naturally extended to the tensor product as 
$$
\tilde{\delta}_{L}:= (\cdot\otimes\id\otimes\id)\circ(\id\otimes\tau\otimes\id)\circ(\delta_{L}\otimes\delta_{L})
$$ 
with $\tau$ the standard flip. A calculus is left covariant if and only if $\tilde{\delta}_{L}(\cn)\subset\ch\otimes\cn$. In this case, the coaction $\delta_{L}^{(1)}$ is the consistent restriction of 
$\tilde{\delta}_{L}$ to $\Omega^{1}(\ca)$.
 The property of right covariance of a first order differential calculus is stated in complete analogy with respect to a right $\ch$-comodule structure of $\ca$.

Clearly these notions make sense for $\ca$ to be the algebra $\ch$ with the coaction $\Delta$ of $\ch$ on itself extended 
then to maps 
\begin{align*}
\Delta^{(1)}_{R}(\dd h)=(\dd\otimes 1)\Delta(h) \qquad \textup{and} \qquad
\Delta^{(1)}_{L}(\dd h)=(1\otimes \dd)\Delta(h).
\end{align*}
On $\ch$ there is in addition the notion of a bicovariant calculus, namely a calculus which is both left and right covariant 
and satisfying the  compatibility condition:
$$
(\id\otimes\Delta^{(1)}_{R})\circ\Delta^{(1)}_{L}=(\Delta^{(1)}_{L}\otimes\id)\circ\Delta^{(1)}_{R}.
$$

On a quantum group $\ch$ the covariance of calculi are studied in terms of the  bijection $r : \ch\otimes\ch\to\ch\otimes\ch,$
\beq
r(h\otimes h^{\prime})=(h\otimes1)\Delta(h^{\prime}) ,\quad\qquad r^{-1}(h\otimes h^{\prime})=(h\otimes1)(S\otimes\id)\Delta(h^{\prime}) 
\label{3p2}
\eeq
which is such that $ r(\och{1}_{un})=\ch\otimes \ker\eps $. 
Left covariant calculi on $\ch$ are  in one to one correspondence with  right ideals 
$\mathcal{Q} \subset \ker\eps$, with subbimodule 
$\mathcal{N}_{\mathcal{Q}}=r^{-1}(\ch\otimes\mathcal{Q})$ and  
$\och{1}:=\och{1}_{un}/\mathcal{N}_{\mathcal{Q}}$.
The left $\ch$-modules isomorphism given by 
$\Omega^{1}(\ch)\simeq\ch\otimes(\ker\eps/\cq)$ allows one  to recover the complex vector space $\ker\eps/\cq$ as the set of left invariant 1-forms, namely the elements $\omega_{a}\in \Omega^{1}(\ch)$ such that 
$$
\Delta^{(1)}_{L}(\omega_{a})=1\otimes\omega_{a} .
$$
The dimension of $\ker\eps/\cq$ is referred to as the dimension of the calculus. 
A left covariant first order differential calculus is  a $*$-calculus if and only if $(S(Q))^*\in \cq$ for any $Q\in \cq$. If this is the case,  the left coaction of $\ch$ on $\Omega^{1}(\ch)$ is compatible with the $*$-structure: 
$\Delta^{(1)}_{L}(\dd h^{*})=(\Delta^{(1)}(\dd h))^*$.
Bicovariant calculi  corresponds to   right ideals $\mathcal{Q}\subset\ker\varepsilon$ which are in addition stable under the right adjoint coaction $\Ad$ of $\ch$ onto itself, that is to say $\Ad(\cq) \subset \cq \otimes \ch$.  
Explicitly, $\Ad=\left(\id \otimes m \right) \left(\tau\otimes
\id \right)\left(S\otimes\Delta\right)\Delta$, with $\tau$ the flip
operator, or $\Ad(h) = \co{h}{2} \otimes \left(S(\co{h}{1}) \co{h}{3} \right)$
in Sweedler notation.

\bigskip
The tangent space of
the calculus is  the complex vector space of elements  out of $\ch^{\prime}$ -- the dual space $\ch^{\prime}$ of functionals on $\ch$ --  defined by 
\beq
\label{qTv}
\mathcal{X}_{\mathcal{Q}}:=
\{X\in \ch^{\prime} ~:~ X(1)=0,\,X(Q)=0, \,\, \forall \, Q\in\mathcal{Q}\}. 
\eeq
There exists a unique bilinear form 
\beq
\hp{~}{~}:\mathcal{X}_{\cq}\times\Omega^{1}(\ch), \qquad\hp{X}{x \dd y}:=\eps(x)X(y) ,
\label{dptg}
\eeq
giving a non-degenerate dual pairing between the vector spaces $\mathcal{X}_{\cq}$ and $\ker\eps/\cq$.  
We have then also a vector space isomorphism $\mathcal{X}_{\cq}\simeq(\ker\varepsilon/\cq)$. 

The dual space $\ch^{\prime}$ has natural left and right (mutually commuting) actions  on $\ch$:  
\beq
X\triangleright h:=h_{(1)}X(h_{(2)}),\qquad
 h\triangleleft X:=X(h_{(1)})h_{(2)}.
\label{deflr}
\eeq
If the  vector space  $\mathcal{X}_{\cq}$ is finite dimensional  its elements belong to the dual Hopf algebra 
$ \ch^{\prime} \supset \ch^{o} = (\ch^{o}, \Delta_{\ch^o}, \eps_{\ch^o},S_{\ch^o})$, defined as the largest Hopf $*$-subalgebra contained in $\ch^{\prime}$. 
In such a case  the 
$*$-structures are compatible with both actions:
$$
X \lt h^* = ((S(X))^* \lt h)^*,\qquad
h^* \rt  X = (h \rt  (S(X))^*)^*,
$$
for any $ X \in \ch^{o}, \ h \in \ch$. 
Then the exterior derivative can be written as:
\beq
\dd h := \sum_a
~(X_{a} \triangleright h) ~\omega_{a} =\sum_{a} \omega_{a}(-S^{-1}(X_{a}))\lt h, 
\label{ded}
\eeq
where $\hp{X_{a}}{\omega_{b}}=\delta_{ab}$, and one has the identity $S^{-1}(X_{a})=-S^{-1}(f_{ba})X_{b}$.  
The twisted Leibniz rule of derivations of the basis elements $X_{a}$ is dictated by their coproduct:  
\beq
\Delta_{\ch^o}(X_{a})=1\otimes X_{a}+
\sum\nolimits_b X_{b}\otimes f_{ba},
\label{cpuh}
\eeq
where the $f_{ab} \in \ch^{o}$ consitute an algebra representation of $\ch$:
\begin{align*}
&\Delta_{\ch^{o}}(f_{ab})=\sum_c f_{ac}\otimes f_{cb}, \nn \\
& \varepsilon_{\ch^{o}}(f_{ab})=\delta_{ab}, \nn \\
&\sum_b S_{\ch^{o}}(f_{ab})f_{bc}=\sum_b f_{ab}S_{\ch^{o}}(f_{bc})=\delta_{ac}.
\end{align*}
 The elements $f_{ab}$ also control the $\ch$-bimodule structure of $\Omega^{1}(\ch)$:
\begin{align}\label{bi-struct}
\omega_{a} h = \sum\nolimits_b (f_{ab}\triangleright h)\omega_{b} , \qquad 
h \omega_{a} = \sum\nolimits_b \omega_{b} \left( (S^{-1}(f_{ab}) )\triangleright h \right) , 
\qquad \mathrm{for} \quad h \in \ch. 
\end{align}

The right coaction of $\ch$ on  $\Omega^{1}(\ch)$ defines matrix elements 
\beq \label{ri-co-form}
\Delta^{(1)}_{R}(\omega_{a})=\sum_b \omega_{b}\otimes J_{ba}.
\eeq
where $J_{ab}\in\ch$. This matrix is invertible, since $\sum_b S(J_{ab})J_{bc}=\delta_{ac}$ and 
$\sum_b J_{ab}S(J_{bc})=\delta_{ac}$. In addition one finds that $\Delta(J_{ab})=\sum_c J_{ac}\otimes J_{cb}$ and 
$\varepsilon(J_{ab})=\delta_{ab}$.  It gives a basis of right invariant 1-forms, $\eta_{a}=\omega_{b}S(J_{ba})$ and, 
as we shall see in a moment, allows one for an explicit evaluation of the braiding of the calculus.


\bigskip

In order to construct an exterior algebra $\Omega(\ch)$ over the bicovariant first order differential calculus $(\Omega^{1}(\ch),\dd)$ one uses a  braiding map replacing the flip automorphism. Define the tensor product $\Omega^{1}(\ch)^{\otimes k}=\Omega^{1}(\ch)\otimes_{\ch}\ldots\otimes_{\ch}\Omega^{1}(\ch)$
with $k$ factors. 
There exists a unique $\ch$-bimodule homomorphism 
$\sigma:\Omega^{1}(\ch)^{\otimes2}\to\Omega^{1}(\ch)^{\otimes2}$ 
such that $\sigma(\omega\otimes\eta)=\eta\otimes\omega$ for any left invariant 1-form $\omega$ and any right invariant 1-form $\eta$. The map $\sigma$ is invertible and commutes with the left coaction of $\ch$: 
$$
(\id\otimes\sigma)\circ\Delta_{L}^{(2)}=\Delta_{L}^{(2)}\circ\sigma,
$$ 
with $\Delta_{L}^{(2)}$ the extension of the coaction to the tensor product. There is an analogous invariance for the right 
coaction. 
Moreover, $\sigma$ satisfies a braid equation. On $\Omega^{1}(\ch)^{\otimes3}$: 
$$
(\id\otimes\sigma)\circ(\sigma\otimes \id)\circ(\id\otimes \sigma)=(\sigma\otimes \id)\circ(\id\otimes \sigma)\circ(\sigma\otimes \id).
$$
All of this was proved in \cite{wor89}, where, using the dual pairing between $\ch^{o}$ and $\ch$, an explicit form of the braiding $\sigma$ was given on a basis of left invariant 1-forms:
\beq
\sigma(\omega_{a}\otimes\omega_{b}):=\sum_{n k}\sigma_{ab}^{\,\,\,\,nk}\omega_{n}\otimes\omega_{k} = \sum_{n k} \hs{f_{ak}}{J_{nb}}\omega_{n}\otimes\omega_{k}.
\label{sigco}
\eeq
The braiding map provides a representation of the braid group and an antisymmetrizer operator
$\mathfrak{A}^{(k)}:\Omega^{1}(\ch)^{\otimes k}\to\Omega^{1}(\ch)^{\otimes k}$. The Hopf ideals $\mathcal{S}_{\cq}^{(k)}=\ker\,\mathfrak{A}^{(k)}$ give the quotients
\beq
\Omega^{k}(\ch)=\Omega^{1}(\ch)^{\otimes k}/\mathcal{S}^{(k)}_{\cq}
\label{wedk}
\eeq
the structure of a $\ch$-bicovariant bimodule  which can be written as $\Omega^{k}(\ch)=\mathrm{Range}\,\mathfrak{A}^{(k)}$. The exterior algebra is  
$(\Omega(\ch)=\oplus_{k}\Omega^{k}(\ch),\wedge)$
with the identification $\Omega^{0}(\ch)=\ch$.  The exterior derivative is extended to $\Omega(\ch)$ as the only degree one derivation such that $\dd^2=0$. The  algebra $\Omega(\ch)$ has natural left and right $\ch$-comodule structure, given by recursively setting
$$
\Delta_{L}^{(k+1)}(\dd\theta)=(1\otimes\dd)\Delta_{L}^{(k)}(\theta), \qquad \Delta_{R}^{(k+1)}(\dd\theta)=(\dd\otimes 1)\Delta_{R}^{(k)}(\theta).
$$
Finally, the $*$-structure on $\Omega^{1}(\ch)$ is  extended to an antilinear $*:\Omega(\ch)\to\Omega(\ch)$ by $(\theta\wedge
 \theta^{\prime})^{*}=(-1)^{kk^{\prime}}\theta^{\prime*}\wedge\theta^{*}$ with $\theta\in \Omega^{k}(\ch)$ and $\theta^{\prime}\in \Omega^{k^{\prime}}(\ch)$; the exterior derivative operator satisfies the identity $(\dd\theta)^{*}=\dd(\theta^*)$. 

The quantum tangent space $\mathcal{X}_{\cq}$ can be endowed with a bilinear product, given as the functional $[~,~]_{q}:\mathcal{X}_{\cq}\otimes\mathcal{X}_{\cq}\to\mathcal{X}_{\cq}$:
\beq
\label{q-comm}
[X, Y]_{q}(h):=\{X \otimes Y,\mathrm{Ad}(h)\}, 
\eeq
with a natural extension of the bilinear form \eqref{dptg}. The bicovariance of the calculus ensures that the product is well defined and that, beside being braided antisymmetric it satisfies a braded Jacobi identity, both properties with respect to the (transpose of) the braiding $\sigma$. On a basis it is given by 
\beq
\label{q-comm-bis}
[X_a, X_b]_{q}=X_{a}X_{b}- \sum_{cd}\sigma_{cd}^{\,\,\,\,ab} X_c X_d  ,
\eeq
and computed in terms of the pairing and the matrix 
$J_{ab}$ in \eqref{ri-co-form} as
\beq\label{q-comm-ter}
X_a X_b - \sum_{cd}\sigma_{cd}^{\,\,\,\,ab} X_c X_d = \sum_{c} \hs{X_{b}}{J_{ac}} X_c .
\eeq

\section{Quantum principal bundles and connections on them}
\label{ap:qpb}

Following \cite{bm93}, we consider as a  total space  an algebra $\mathcal{P}$ (with multiplication $m : \mathcal{P}\otimes \mathcal{P} \to \mathcal{P}$) and as structure group  
a Hopf algebra $\ch$, thus $\mathcal{P}$ is a right $\ch$-comodule algebra
with coaction $\delta_{R}\,:\,\mathcal{P}\to\mathcal{P}\otimes\ch$.
The subalgebra of  right coinvariant elements, 
$
\mathcal{B}=\mathcal{P}^{\ch}\,=\,\{p\in\mathcal{P}\,:\,\delta_{R} p
= p\otimes 1\}
$,
is the base space of the bundle. The algebras $(\mathcal{P}, \mathcal{B}, \mathcal{H})$ define a topological principal bundle provided the  
sequence:
\beq
0\, \to \,\mathcal{P}\left(\Omega^{1}(\mathcal{B})_{un}\right)\mathcal{P}\,
 \to \,\Omega^{1}(\mathcal{P})_{un}
\,\stackrel{\chi}  \to \,\mathcal{P}\otimes \ker\varepsilon_{\ch}\,  \to \,0
\label{topes}
\eeq
is exact, with $\Omega^{1}(\mathcal{P})_{un}$ and $\Omega^{1}(\mathcal{B})_{un}$ the universal calculi and the map $\chi$  defined by
\beq
\chi:\,\mathcal{P}\otimes\mathcal{P}\, \to\,\mathcal{P}\otimes\ch,\qquad 
\chi:=\left(m \otimes \id\right)\left(\id\otimes\delta_{R}\right)
\label{chimap}.
\eeq
In fancier parlance, the exactness of this sequence is also referred to as stating (for 
a structure quantum group which is cosemisimple and has bijective antipode)
that the inclusion $\mathcal{B} \hookrightarrow  \mathcal{P}$ is a Hopf-Galois extension \cite[Th.~I]{sc90}. 

Assume now that $\left(\Omega^{1}(\mathcal{P}),\dd\right)$ is 
a right $\ch$-covariant differential calculus  on $\mathcal{P}$ given via the subbimodule $\mathcal{N}_{\mathcal{P}} \subset \Omega^{1}(\mathcal{P})_{un}$, and $\left(\Omega^{1}(\ch),\dd\right)$ a bicovariant differential calculus on $\ch$ given via the $\Ad$-invariant right ideal $\mathcal{Q}_{\ch} \in \ker\varepsilon_{\ch}$. A first order left invariant differential calculus is induced on the algebra basis $\mathcal{B}$ via $\Omega^{1}(\mathcal{B})=\Omega^{1}(\mathcal{B})_{un}/\cn_{\mathcal{B}}$ with 
$\mathcal{N}_{\mathcal{B}}:= \mathcal{N}_{\mathcal{P}} \cap \Omega^{1}(\mathcal{B})_{un}$. This definition is aimed to ensure that  $\Omega^{1}(\mathcal{B})=\cb \dd \cb$. 

To extend the coaction $\delta_{R}$ to a coaction of $\ch$ on $\Omega^{1}(\cp)$, one requires $\delta_{R}(\cn_{\cp})\subset\cn_{\cp}\otimes \ch$.
The compatibility of the calculi are then the  requirements that 
$\chi(\mathcal{N}_{\mathcal{P}})\subseteq \mathcal{P}\otimes\mathcal{Q}_{\ch}$ and that the map 
$\sim_{ {\mathcal{N}_{{\mathcal{P}} }}} : \Omega^{1}(\mathcal{P}) \to  
\mathcal{P}\otimes(\ker\varepsilon_{\ch}/\mathcal{Q}_{\ch})$, defined by the diagram
\beq
\begin{array}{lcl}
\Omega^{1}(\mathcal{P})_{un} 
& \stackrel{{\pi_{\mathcal{N}}}} {\longrightarrow} & \Omega^{1}(\mathcal{P}) \\
\downarrow \chi  &  & \downarrow\sim_{{ \mathcal{N}_{\mathcal{P}}}} \\
\mathcal{P}\otimes \ker\varepsilon_{\ch} &
\stackrel{{\id\otimes\pi_{\cq_{\ch}}}}{\longrightarrow}&
\mathcal{P}\otimes(\ker\varepsilon_{\ch}/\mathcal{Q}_{\ch})
\end{array}
\label{qdia}
\eeq
(with $\pi_{\cn}$ and $\pi_{\cq_{\ch}}$ the natural projections), is
surjective and has kernel 
$\ker(\sim_{{\mathcal{N}_{\mathcal{P}}}}) = \mathcal{P}\Omega^{1}(\mathcal{B})\mathcal{P}
=:\Omega^{1}_{\mathrm{hor}}(\cp)$.
These conditions ensure the exactness of the sequence:
\beq
0\,\to\,\mathcal{P}\Omega^{1}(\mathcal{B})\mathcal{P}\,\to\,
\Omega_{1}(\mathcal{P})
\,\stackrel{\sim_{\mathcal{N}_{\mathcal{P}}}} \longrightarrow\, \mathcal{P}\otimes \left(\ker\varepsilon_{\ch}/\mathcal{Q}_{\ch}\right)\,\to\,0 .
\label{des}
\eeq
The condition 
$\chi(\mathcal{N}_{\mathcal{P}})\subseteq \mathcal{P}\otimes\mathcal{Q}_{\ch}$ is needed to have a well defined map 
$\sim_{ {\mathcal{N}_{{\mathcal{P}} }}}$. If $(\cp, \cb, \ch)$ is a quantum principal bundle with the universal calculi, the equality $\chi(\mathcal{N}_{\mathcal{P}})=\mathcal{P}\otimes\mathcal{Q}_{\ch}$ ensures that $(\cp, \cb, \ch; \cn_{\cp}, \cq_{\ch})$ is a quantum principal bundle with the corresponding nonuniversal calculi.


Elements in the quantum tangent space  $\mathcal{X}_{\cq_{\ch}}(\ch)$ act on $\ker\varepsilon_{\ch} / \mathcal{Q}_{\ch}$
via the pairing  between $\ch^{o}$ and $\ch$. Given  $V\in\mathcal{X}_{\cq_{\ch}}(\ch)$ one defines a map 
\beq\label{hf}
\widetilde{V} : \Omega^{1}(\mathcal{P}) \to\mathcal{P}, \qquad \widetilde{V}:=\left(\id\otimes V\right) \circ (\sim_{\mathcal{N}_{\mathcal{P}}} )
\eeq
and declares a 1-form $\omega \in \Omega^{1}(\mathcal{P})$ to be  horizontal if{}f
$\widetilde{V}\left(\omega\right)=0$,  for any $V\in
\mathcal{X}_{\cq_{\ch}}(\ch)$.
The collection of horizontal 1-forms  coincides with $\Omega^{1}_{\mathrm{hor}}(\cp)$.

\bigskip
The compatibility conditions above allow one to define  right coactions (for $k=0,1,\ldots$) 
$\delta_{R}^{(k+1)}:\Omega^{k+1}(\cp)\to\Omega^{k+1}(\cp)\otimes\ch$, as  coalgebra maps, via $\delta_{R}^{(k+1)}\circ\dd=(\dd\otimes 1)\circ\delta^{(k)}_{R}$.
By direct computation $\Ad(\ker\varepsilon_{\ch})\subset(\ker\varepsilon_{\ch})\otimes\ch$. Being  the right ideal $\cq_{\ch}$  $\Ad$-invariant (i.e. the differential calculus on $\ch$ is bicovariant),  it is possible to define a right-adjoint coaction $\Ad^{(R)}: \ker\varepsilon_{\ch}/\cq_{\ch} \,\to\,(\ker\varepsilon_{\ch}/\cq_{\ch})\otimes\ch$ by the commutative diagram
$$
\begin{array}{lcl}
\ker\varepsilon_{\ch} 
& \stackrel{{\pi_{\cq_{\ch}}}} {\longrightarrow} & \ker\varepsilon_{\ch}/\cq_{\ch} \\
\downarrow \Ad  &  & \downarrow \Ad^{(R)}   \\
\ker\varepsilon_{\ch}\otimes\ch &
\stackrel{{\pi_{\cq_{\ch}}}\otimes\id}{\longrightarrow}&
(\ker\varepsilon_{\ch}/\mathcal{Q}_{\ch})\otimes\ch
\end{array}
$$
Such a right-adjoint coaction $\Ad^{(R)}$ allows one further to define a right coaction $\delta_{R}^{(\Ad)}$ of $\ch$ on $\cp\otimes(\ker\varepsilon_{\ch}/\cq_{\ch})$ as a coaction of a Hopf algebra on the tensor product of its comodules. This coaction is explicitly given by the relation:
\beq
\delta_{R}^{(\Ad)}(p\otimes\pi_{\cq_{\ch}}(h))=p_{(0)}\otimes\pi_{\cq_{\ch}}(h_{(2)})\otimes p_{(1)}(Sh_{(1)})h_{(3)}.
\label{Adre}
\eeq

A connection  
on the quantum principal bundle is a right invariant splitting of the sequence \eqref{des}. Given a left $\cp$-linear map $\sigmat:\cp\otimes(\ker\varepsilon_{\ch}/\cq_{\ch})\to\Omega^{1}(\cp)$ such that
\begin{align}
 \delta_{R}^{(1)}\circ\sigmat=(\sigmat\otimes \id)\delta^{(\Ad)}_{R} ,\qquad \textup{and} \qquad
 \sim_{\cn_{\cp}}\circ\sigmat=\id, 
\label{si}
\end{align}
the map $\Pi:\Omega^{1}(\cp)\to\Omega^{1}(\cp)$ defined by $\Pi=\sigmat\ \circ\sim_{\cn_{\cp}}$ is a right invariant left $\cp$-linear projection, whose kernel coincides with the horizontal forms $\mathcal{P}\Omega^{1}(\mathcal{B})\mathcal{P}$: 
\begin{align}
&\Pi^{2}=\Pi ,\nn \\
&\Pi(\mathcal{P}\Omega^{1}(\mathcal{B})\mathcal{P})=0,\nn\\
&\delta_{R}^{(1)}\circ\Pi=(\Pi\otimes \id)\circ\delta_{R}^{(1)}.
\label{Pi}
\end{align}
The image of the projection $\Pi$ is the set of vertical 1-forms of the principal bundle. A connection on a principal bundle can also be given via a connection one form, which is a map $\omega:\ch\to\Omega^{1}(\cp)$. Given a right invariant splitting $\sigmat$ of the exact sequence \eqref{des},  define the connection 1-form as $\omega(h)=\sigmat(1\otimes\pi_{\cq_{\ch}}(h-\varepsilon_{\ch}(h)))$ on $h\in\,\ch$. Such a  connection 1-form has the following properties:
\begin{align}
&\omega(\cq_{\ch})=0, \nn \\
&\sim_{\cn_{\cp}}(\omega(h))=1\otimes\pi_{\cq_{\ch}}(h-\varepsilon_{\ch}(h))\qquad\forall\,h\in\,\ch,\nn \\
&\delta_{R}^{(1)}\circ\omega=(\omega\otimes\id)\circ \Ad,\nn \\
&\Pi(\dd p)=\cdot\, (\id\otimes\omega)\delta_{R}(p)\qquad\forall\, p\in\,\cp.
\label{ome} 
\end{align}
 Conversely with a linear map $\omega:\ker\varepsilon_{\ch}\to\Omega^{1}(\cp)$ that satisfies the first three conditions in \eqref{ome},  there exists a unique connection on the principal bundle, such that $\omega$ is its connection 1-form. The splitting of the sequence \eqref{des} will be 
\beq
\sigmat(p\otimes[h])=p\omega([h])
\label{siom}
\eeq
with $[h]$ in $\ker(\varepsilon_{\ch}/\cq_{\ch})$, while the projection $\Pi$ will be 
\beq
\Pi=m\circ(\id\otimes\omega)\circ\sim_{\cn_{\cp}}
\label{piom}
\eeq
The general proof of these results is in \cite{bm93}.

\end{document}